\documentclass[11pt,a4paper]{article}
\usepackage[T1]{fontenc}
\usepackage{lmodern}
\usepackage[margin=29mm]{geometry}
\usepackage{amsmath,amssymb,amsthm,mathtools,booktabs,longtable,array,enumitem,microtype,tikz}
\usepackage{needspace}
\usepackage[hidelinks]{hyperref}
\usepackage[font=small,labelfont=bf]{caption}
\hypersetup{pdftitle={Positive recurrence of weakly reversible stochastic reaction networks with mass-action kinetics},pdfauthor={Daniele Cappelletti and Giulio Cuniberti}}
\newcommand{\nat}{\mathbb Z_{\ge0}}
\newcommand{\pathcost}{\operatorname{cost}}
\newcommand{\Sp}{\mathcal{S}}
\providecommand{\C}{}
\renewcommand{\C}{\mathcal{C}}
\renewcommand{\Re}{\mathcal{R}}
\renewcommand{\L}{\mathcal{L}}
\newcommand{\dtp}[1]{\langle #1\rangle}
\newtheorem{theorem}{Theorem}[section]
\newtheorem{lemma}[theorem]{Lemma}
\newtheorem{proposition}[theorem]{Proposition}
\theoremstyle{definition}
\newtheorem{definition}[theorem]{Definition}
\newtheorem{example}[theorem]{Example}
\newtheorem{remark}[theorem]{Remark}
\newtheorem*{examplecont}{Example~\ref{ex:running} (continued)}
\theoremstyle{plain}
\setlist{topsep=0.4em,itemsep=0.25em,parsep=0pt}
\numberwithin{equation}{section}
\allowdisplaybreaks[1]
\title{\bfseries Positive recurrence of weakly reversible\\stochastic reaction networks\\with mass-action kinetics}
\author{Daniele Cappelletti and Giulio Cuniberti\\[0.5em]{\small Department of Mathematical Sciences, Politecnico di Torino, Turin, Italy}\\[0.2em]{\small \href{mailto:daniele.cappelletti@polito.it}{\texttt{daniele.cappelletti@polito.it}}\quad \href{mailto:giulio.cuniberti@polito.it}{\texttt{giulio.cuniberti@polito.it}}}}
\date{26 November 2026}
\begin{document}
\maketitle

\begin{abstract}
We prove that for every finite weakly reversible stochastic reaction network with mass-action kinetics and arbitrary positive rate constants, the Markov chain is non-explosive and positive recurrent on each communicating class. This settles the long-standing chemical recurrence conjecture. Weak reversibility means that every reaction lies on a directed cycle of complexes. The proof constructs a Lyapunov function explicitly. We first fix nonnegative weights on the species. To each path of reactions we assign the largest among the weighted counts of the molecules left unused at its steps and the weighted count of its final state plus a fixed positive constant. The height of a state is the minimum of this quantity over the paths starting from that state. A reaction that lowers the height may become active only after other reactions have fired. We add to the height a bounded correction that measures how many of these preparatory reactions are still needed, and we use weak reversibility to control the reactions that undo part of the preparation. Finally, we sum the resulting functions over finitely many choices of the weights, and we show that the sum has negative drift outside a finite set. We prove the result more generally for variable-$\kappa$ mass-action kinetics, in which the rate constants are replaced by functions of the state that are bounded from above and away from zero.
\end{abstract}



\paragraph{Contributions and use of AI.} The authors directed the research, developed and refined parts of the proof, checked the mathematical arguments, and revised the manuscript. They used GPT-6 Astra (OpenAI) to develop mathematical constructions and arguments, investigate possible counterexamples, perform internal checks, and draft the text. Claude Opus 5.5 and Fable 5.2 (Anthropic), together with GPT-6.2 Astra (OpenAI), helped review drafts for mathematical correctness and clarity. The authors are responsible for the manuscript.

\section{Introduction}\label{sec:introduction}

Stochastic reaction networks describe how the counts of interacting species change over time~\cite{AndersonKurtz2015,Toth2018}. The state of the model is the vector of the molecule counts, and it evolves as a continuous-time Markov chain whose jumps correspond to the reactions. In this paper we study when such a chain is positive recurrent, that is, when it returns to its starting state in finite mean time after leaving it. A closed positive recurrent communicating class carries a unique stationary distribution, and if the chain starts in that class, then the fraction of time spent in each state converges almost surely to its stationary probability~\cite{Norris1997}.

As an illustration, consider the network $0\rightleftharpoons A$. Molecules of $A$ are produced at a constant rate, and each of them is removed at a constant rate. Hence, the total removal rate grows with the count, and the chain is pushed back whenever the count is large. With several species, this mechanism can fail at the boundary of the state space: a reaction that would lower a large count may need a molecule of another species whose count is zero, and in this case it cannot fire until that molecule is produced. This is the main difficulty in proving positive recurrence for general networks.

The \emph{chemical recurrence conjecture}, also called the \emph{positive recurrence conjecture}, states that, under mass-action kinetics, the chain associated with a weakly reversible network is positive recurrent on each communicating class, for every choice of positive rate constants. A network is weakly reversible if every reaction lies on a directed cycle of reactions (Definition~\ref{def:network}). Anderson and Kim state the conjecture in~\cite{AK2018}. Complex-balanced systems admit product-form stationary distributions~\cite{ACK2010} and are non-explosive~\cite{ACKK2018}. Anderson, Cappelletti, and Kim~\cite{ACK2020} prove positive recurrence for networks with one weakly reversible linkage class, at most two molecules per complex, and a positive multiple of each species occurring as a complex. Kriebel proposes a proof without the latter condition~\cite{Kriebel2026}. Wiuf and Xu cover networks with one-dimensional stoichiometric subspace~\cite{WX2023}. Agazzi, Anderson, Cappelletti, Laurence, and Mattingly give a proof for two species~\cite{Agazzi2025}. Other related results concern planar bulk estimates with separate boundary conditions~\cite{AgazziLaurence2026}, non-explosion~\cite{Xu2026regularity}, and exponential ergodicity for first-order endotactic systems~\cite{Xu2026firstorder}. In this paper we prove the conjecture: for every finite weakly reversible network and for any positive rate constants, the mass-action chain is non-explosive from every initial state and positive recurrent on each communicating class (Theorem~\ref{thm:recurrence}). We state and prove the result for the more general variable-$\kappa$ mass-action kinetics, in which the rate constants are replaced by functions of the state that are bounded from above and away from zero. Mass-action kinetics is the special case of constant functions. Rate constants that vary in time within such bounds define the variable-$\kappa$ mass-action systems studied for the deterministic model~\cite{CNP2013}.

For the deterministic model, the corresponding conjecture asks whether each positive stoichiometric compatibility class contains a compact set in which every trajectory starting in that class eventually remains~\cite{CNP2013}. A proof was recently proposed in~\cite{OpenAI2026}. In the stochastic model the situation is different, because the count of a species can reach zero in a single reaction, and the reactions that need that species then stop. As in the works on tiers and on Lyapunov functions for large molecule counts~\cite{ACKN2020,ADE2018,ADEgeometry2018}, our argument compares the orders of magnitude of the reaction rates when the molecule counts are large. At the boundary a further difficulty arises. A decrease of the large counts may require a sequence of reactions that first produce the missing molecules, and a faster reaction can fire in the meantime and undo part of this sequence. A large part of the proof is devoted to controlling this competition.

The proof for two species in~\cite{Agazzi2025} patches together functions defined on different regions of the state space. Here, the Lyapunov function is given by a single formula, and the drift estimate holds at every state outside a finite set.

The paper is organized as follows. Section~\ref{sec:model} collects the definitions, and Section~\ref{sec:main} states the main result and outlines the proof. In Section~\ref{sec:height} we fix weights on the species and define the height of a state by means of the paths of reactions starting from it. A path is allowed to stop at any state, and it is charged for the molecules present at its endpoint. With this choice, when more molecules are added to every state of a path, the resulting path is still admissible and its charge increases by the weighted count of the added molecules, a property used throughout the proof. In Section~\ref{sec:races} we assign a cost to the reactions that are needed before the height can decrease, and we use the total cost to correct the height. In Section~\ref{sec:selection} we select finitely many vectors of species weights, adapting the affine approximation of~\cite{OpenAIBP2026}. The function associated with each of these vectors is expressed in terms of its own height, and for this reason the functions can be added directly in Section~\ref{sec:assembly}, where we conclude the proof of the main result.

\section{Preliminaries}\label{sec:model}

We write $\nat=\{0,1,2,\ldots\}$. Inequalities and logarithms of vectors are taken coordinatewise. We use the norms $\|v\|_1=\sum_i|v_i|$ and $\|v\|_\infty=\max_i|v_i|$, and write $\mathbf1$ for a vector of ones. For a vector $v$, we write $v+1$ for the vector with entries $v_i+1$.

\begin{definition}[Reaction network and weak reversibility]\label{def:network}
Fix a set $\Sp=\{S_1,\ldots,S_d\}$ of $d$ \emph{species}, in this order. A \emph{state} $x\in\nat^d$ records how many molecules of each species are present. A \emph{complex} $y\in\nat^d$ describes a combination $y_1 S_1 + \ldots + y_d S_d$ of these same species, with $y_i$ molecules of species $i$. A \emph{reaction} $y\to y'$ consumes the \emph{source complex} $y$ and produces the \emph{product complex} $y'$. A \emph{reaction network} $(\Sp,\C,\Re)$ consists of the species set $\Sp$, a finite set $\C$ of complexes, and a finite set $\Re$ of reactions. Its \emph{reaction graph} has the complexes as vertices and the reactions as directed edges. A \emph{linkage class} is a connected component of the underlying undirected graph. The network is \emph{weakly reversible} if each linkage class is strongly connected, or equivalently if each reaction lies on a directed cycle. We discard reactions with $y=y'$, which leave the state unchanged.
\end{definition}

\begin{definition}[Variable-$\kappa$ mass-action kinetics]\label{def:kinetics}
A reaction $y\to y'$ is \emph{active} at $x$ if $y\le x$, so that all its reactants are present. When it fires, the state becomes $x-y+y'$. Under \emph{stochastic variable-$\kappa$ mass-action kinetics}, each reaction $y\to y'$ has a \emph{rate coefficient} $\kappa_{y\to y'}:\nat^d\to[\kappa_{\min},\kappa_{\max}]$, where the bounds $0<\kappa_{\min}\le\kappa_{\max}<\infty$ do not depend on the reaction, and its rate at $x$ is
\begin{equation}\label{eq:mass-action}
 \lambda_{y\to y'}(x)=\kappa_{y\to y'}(x)\lambda_y(x),\qquad \lambda_y(x)=
 \begin{cases}
     \displaystyle\prod_{i=1}^d\frac{x_i!}{(x_i-y_i)!} & \text{if }y\le x, \\[1em]
     0 & \text{otherwise}.
 \end{cases}
\end{equation}
The factor $x_i!/(x_i-y_i)!=x_i(x_i-1)\cdots(x_i-y_i+1)$ counts ordered choices of $y_i$ molecules of species $S_i$ from the $x_i$ available.
\end{definition}

\begin{remark}[Mass-action kinetics]\label{rem:mass-action}
The usual \emph{stochastic mass-action kinetics} is the special case in which every rate coefficient is a positive constant $\kappa_{y\to y'}$, called the \emph{rate constant} of the reaction. In this case, one can take $\kappa_{\min}$ and $\kappa_{\max}$ to be the smallest and the largest rate constant.
\end{remark}

For example, with species ordered as $A,B$, the source complex $2A+B$ is the vector $y=(2,1)$. A reaction $y\to y'$ from this complex has rate $\kappa_{y\to y'}(x) x_A(x_A-1)x_B$. To be active, it needs two molecules of $A$ and one of $B$.

We follow the successive jumps of the \emph{minimal continuous-time Markov chain} $X$ up to their possible accumulation time. Part of the proof is to show that this time is almost surely infinite. Let $\L $ be its \emph{generator}. For $f:\nat^d\to\mathbb R$,
\begin{equation}\label{eq:generator}
 \L f(x)=\sum_{\substack{y\to y'\in\Re\\ y\le x}}\lambda_{y\to y'}(x)\bigl[f(x-y+y')-f(x)\bigr].
\end{equation}
Thus, $\L f(x)$ is the sum, over all active reactions, of the rate times the change in $f$. We call this the \emph{drift} of $f$ at $x$.

\begin{definition}[Communicating classes and recurrence]\label{def:classes}
Two states \emph{communicate} if each can be reached from the other through finitely many active reactions. The equivalence classes for this relation are the \emph{communicating classes}. Such a class is \emph{closed} if no active reaction can take the chain out of it. The chain is \emph{irreducible} if its state space forms a single communicating class. It is called \emph{non-explosive} if it has finitely many jumps on every bounded time interval, almost surely. A state that has at least one active reaction is \emph{positive recurrent} if the expected time to return to it after the first departure is finite. A return that never occurs has time $\infty$. A state with no active reaction is \emph{absorbing}. The chain stays there forever, and we count the state as positive recurrent by convention. Positive recurrence is a class property, meaning that if one state is positive recurrent, then so is every other state in its communicating class~\cite{Norris1997}.
\end{definition}

\begin{example}[Constant production and linear removal]\label{ex:warmup}
For $0\rightleftharpoons A$, with constant rate coefficients, production rate $\kappa_1$, and removal rate $\kappa_2x$, take $V(x)=e^x$. Then
\[
 \frac{\L V(x)}{V(x)}=\kappa_1(e-1)+\kappa_2x(e^{-1}-1)\longrightarrow-\infty.
\]
The negative removal term eventually dominates production. The proof below recovers this effect for interacting species by correcting an exponential of a weighted count.
\end{example}

\section{Main result and outline of the proof}\label{sec:main}

\begin{theorem}\label{thm:recurrence}
Consider a finite weakly reversible reaction network with variable-$\kappa$ mass-action kinetics. Then the minimal chain $X$ is non-explosive from every initial state. Every communicating class is closed, and the chain restricted to each such class is positive recurrent.
\end{theorem}

By Remark~\ref{rem:mass-action}, Theorem~\ref{thm:recurrence} covers mass-action kinetics with arbitrary positive rate constants. We fix the network and its rate coefficients throughout. The estimates use the rate coefficients only through the bounds $\kappa_{\min}$ and $\kappa_{\max}$. Constants in our estimates may depend on these choices but never on the varying state. Appendix~\ref{app:notation} summarizes the other dependencies.

Suppose that a reaction $y\to y'$ is active at a state $x$, so that $x-y\ge0$ collects the molecules that are not consumed by it. The reaction takes $x$ to $x-y+y'$. Weak reversibility gives a directed path in the reaction graph from $y'$ back to $y$. At each step leave the same $x-y$ molecules aside. The remaining molecules initially have counts $y'$, and each reaction along the path produces exactly the reactants needed for the next. At the end of the path the state is $(x-y)+y=x$. Every active transition thus admits a feasible path back to its source state, so every communicating class is closed~\cite{PCK2014}. Restricted to a finite communicating class, the chain is finite and irreducible, hence non-explosive and positive recurrent~\cite{Norris1997}. From now on, we fix an infinite communicating class $\Gamma$.

Our target is a function $V:\Gamma\to(0,\infty)$ with \emph{finite sublevel sets}, meaning that $\{x\in\Gamma:V(x)\le \alpha\}$ is finite for every finite $\alpha$, and for which
\begin{equation}\label{eq:target-drift}
 \L V(x)\le-1\qquad\text{outside a finite subset of }\Gamma .
\end{equation}
Condition \eqref{eq:target-drift} states that, outside a finite set, the expected rate of change of $V$ along the chain is at most $-1$. Finite sublevel sets ensure that $V$ tends to infinity along every sequence of distinct states. By the Foster criterion (Proposition~\ref{prop:foster}), the existence of such a function implies that $X$ is positive recurrent on $\Gamma$. At the end of the proof we show that the repeated returns also rule out explosion. Any such $V$ is a \emph{Lyapunov function} for $X$.

In the example $0\rightleftharpoons A$, we could see why the chain comes back by looking only at the number of molecules of $A$. With several species there is a new possibility. A reaction that would bring a large count down may require one molecule of another species whose count is zero. Producing that molecule can be a step toward a return, even if the large count has yet to fall. We want a function $V$ that accounts for both kinds of progress.

\paragraph{Weighted counts.} The function $V$ is built from weighted molecule counts $\dtp{w,x}=\sum_i w_ix_i$, where the vector $w\in[0,1]^d$ assigns a weight to each species and is called a \emph{direction}. Which weights are useful at a state $x$ depends on how the counts of the species compare with each other, and we measure this on a logarithmic scale. Let $T=\sum_i(x_i+1)$. Adding one to each count lets us take logarithms even when a count is zero. For $T>1$, we call $\log_T (x_i+1)$ the \emph{exponent} of species $i$ at $x$, so that $x_i+1$ is $T$ raised to this exponent. For example, if the counts are $(n,\lfloor\sqrt n\rfloor,3)$ and $n$ grows, then $T/n\to1$ and the exponents tend to $1,\tfrac12,0$. In this case, the direction that describes these states is $w=(1,\tfrac12,0)$. The reason for matching the weights with the exponents is that, when $x_i+1$ is close to $T^{w_i}$ for every $i$, the rate of an active reaction with source complex $y$ is of order $T^{\dtp{w,y}}$, so $w$ determines which reactions are the fastest. The curve $\theta\mapsto(\theta^{w_1},\ldots,\theta^{w_d})$ is the \emph{toric ray} with direction $w$ of~\cite{GMS2014}, a special case of the \emph{toric jets} used for large deviations in~\cite{ADE2018}. The same directions appear in the analysis of the two-species case. In~\cite{AgazziLaurence2026} the scalings $x_i\sim N^{w_i}$ are considered with the same normalization $\max_i w_i=1$. The directions associated with the faces of the convex hull of the source complexes split the large states into strips around these directions and regions in between, and a Poincar\'e compactification is carried out for each of these directions. In~\cite{Agazzi2025} the states $(N,N^{w_2/w_1})$ are considered for $w$ on the unit circle. At the state $x=(10{,}000,100,3)$ the weighted count is $\dtp{w,x}=10{,}000+\tfrac12(100)=10{,}050$. A limiting exponent of zero covers both a bounded count and one that grows more slowly than every positive power of $n$, such as $\lfloor\log n\rfloor$. We work with the exact counts, and not only with their exponents, because even one molecule can make a reaction active.

A single direction cannot describe all the states of $\Gamma$, because the exponents change from state to state. We therefore fix a finite family $\mathcal A$ of directions once for the whole class, and at each state we rely on the directions of the family that are close to the exponents of that state. The family, and the rule that singles out these directions, are described at the end of this section.

\paragraph{Height and progress before a decrease.} Fix a direction $w$ and consider a path of reactions starting from $x$. At each step, the molecules of the current state that are not consumed by the reaction are called \emph{unused}, and we consider their weighted count. At the endpoint $z$ of the path we consider the weighted count $\dtp{w,z}$ plus a fixed positive constant. The \emph{ceiling} of the path is the largest of these numbers. The \emph{height} $H_w(x)$ of the state $x$ is defined as the smallest ceiling among all the paths starting from $x$ (Definition~\ref{def:height}). The \emph{state correction} is the difference between the weighted count and the height,
\[
 s_w(x)=\dtp{w,x}-H_w(x).
\]
Theorem~\ref{thm:height} shows that $s_w$ is bounded. Hence, the height is a bounded modification of the weighted count. When the height is above a fixed threshold, there is a path along which the height stays constant until its last reaction, which lowers it. We call the reactions before the last one \emph{preparatory}. They do not change the height, so the height alone does not register that they bring the chain closer to a decrease. To this aim, we add to the height a bounded correction that measures the preparation still needed.

While the preparation is in progress, a different reaction can fire and lead the chain away from the path. In this case, weak reversibility supplies a path back to the state occupied before that reaction, along which the unused molecules are the same at every step. Hence, the additional preparation caused by such a detour is bounded, even if the detour is fast, as we show in Section~\ref{sec:races}.

\paragraph{Choice of the directions.} It remains to explain how the directions of the family $\mathcal A$ are used at a given state. We say that a direction $w$ \emph{fits} at $x$ if every coefficient $w_i$ is within a prescribed tolerance of the exponent $\log_T (x_i+1)$. The drift estimate of Section~\ref{sec:races} holds for $w$ at the large states where $w$ fits. To single out the fitting directions, we subtract from the height of each direction a \emph{penalty}, which charges every species, beyond the first, that receives a positive weight. A positive weight increases the height by an amount proportional to the count of the species, so it compensates for the penalty only if the species is sufficiently abundant. The height minus the penalty is the \emph{score} of the direction. In Section~\ref{sec:selection} we choose the family and the penalties so that, at every large state, every direction whose score is close to the largest one fits. The function $V$ is a sum over $\mathcal A$ of terms proportional to the exponential of the score. Hence, at a large state the terms of the directions with a score close to the largest one have negative drift, and the terms of the remaining directions are too small to offset it.

In full, the function constructed in Section~\ref{sec:assembly} is
\[
 V(x)=\sum_{w\in\mathcal A}e^{H_w(x)-b_w(H_w(x))}\bigl(1+D_w(x)\bigr).
\]
Here $H_w$ is the height associated with the direction $w$ (Section~\ref{sec:height}), $D_w$ is the bounded correction for the remaining preparation (Section~\ref{sec:races}), and $\mathcal A$ is the finite family of selected directions with their penalties $b_w$ (Section~\ref{sec:selection}). The construction thus has four stages: a height, a correction for preparatory steps, the selection of finitely many directions, and their sum (Section~\ref{sec:assembly}).

\section{The height of a state}\label{sec:height}

\subsection{Unused molecules and a terminal charge}
\begin{definition}[Direction and complex order]\label{def:slope}
Fix a direction $w\in[0,1]^d$ with $\max_i w_i=1$. The \emph{weighted count} at state $x$ is $\dtp{w,x}$, and the \emph{complex order} of $y$ is $\dtp{w,y}$.
\end{definition}

The class $\Gamma$ and the direction $w$ are fixed throughout Sections~\ref{sec:height}--\ref{sec:races}. Constants may depend on the network, on the bounds $\kappa_{\min}$ and $\kappa_{\max}$, on $\Gamma$, and on $w$, but not on the state. The quantities introduced in these sections that depend on $w$ carry the subscript $w$, and Appendix~\ref{app:notation} lists the order in which the constants are chosen.

The \emph{nominal rate} of a reaction with source $y$ is $T^{\dtp{w,y}}$, where $T=\sum_i(x_i+1)$. When each $x_i+1$ is close to $T^{w_i}$, the rate of an active reaction is close to its nominal rate, and Lemma~\ref{lem:actual-rates} gives the precise comparison. The exact counts still determine which reactions are active.

The following heuristic explains why the weighted count is expected to decrease in the simplest situation. Suppose that every reaction is active, that the exponents of the state tend to $w$, and that in some linkage class the complex orders are not all equal. In that linkage class, follow a directed path from a complex of highest order to one of lower order, which exists by weak reversibility. The first reaction on this path that leaves the highest order lowers the weighted count. Every reaction of that linkage class that raises the weighted count starts from a complex of lower order, so the ratio between the combined rate of these reactions and the rate of this one decrease tends to zero. Both assumptions of the heuristic can fail. At the boundary of the state space the decrease itself can be inactive, and the complex orders of a linkage class can all be equal, in which case no reaction of that class changes the weighted count. For these reasons, we do not rely on the complex orders to find a decrease. We look instead for a sequence of active reactions that leads to a state with a smaller weighted count, and Lemma~\ref{lem:terminal-choice} shows that such a sequence exists from every state of $\Gamma$ whose weighted count is large enough.

\begin{definition}[Reaction steps and unused molecules]\label{def:step}
A \emph{step} $(x,y\to y')$ consists of a state $x$ and a reaction $y\to y'$ active there. Its \emph{destination} is $x-y+y'$. A \emph{path} from $x$ to $z$ is a finite sequence of steps such that the first one starts at $x$, each of the others starts at the destination of the preceding one, and the last one has destination $z$. The empty path goes from $x$ to $x$. We call $x-y\ge0$ the \emph{unused stock} of the step $(x,y\to y')$. These are the counts of molecules set aside while the reactants $y$ are replaced by the products $y'$. The \emph{unused level} of the step is the weighted count of its unused stock,
\[
 \dtp{w,x-y}=\dtp{w,x}-\dtp{w,y}.
\]
\end{definition}
Weak reversibility provides a directed path of complexes from $y'$ back to $y$. Adding the unused stock $x-y$ to every complex turns it into a path from $x-y+y'$ back to $x$. All its steps have the same unused level $\dtp{w,x-y}$. This is the reason why we measure unused levels: they do not change along the return paths provided by weak reversibility.

\begin{example}[A boundary network]\label{ex:running}
Our running example is
\begin{equation}\label{eq:example-network}
 2A+B\rightleftharpoons A,\qquad B\rightleftharpoons0.
\end{equation}
Write $(n,k)$ for the counts of $A,B$ and take $w=(1,0)$. Adding the stock $(m,k)$ to the complexes of the first linkage class gives
\[
 \underbrace{(m+1,k)}_{A+(m,k)}
 \quad\rightleftharpoons\quad
 \underbrace{(m+2,k+1)}_{2A+B+(m,k)}.
\]
Both steps have unused level $m$. The complex orders $1$ and $2$ belong to $A$ and $2A+B$, respectively, and do not change when the stock is added. At $(n,0)$ with $n\ge1$, producing a $B$ by $A\to2A+B$ leads to $(n+1,1)$, whereas using $0\to B$ leads to $(n,1)$. Both reactions produce a molecule of $B$, but they lead to different states.
\end{example}

We use paths to detect that the weighted count can decrease. A path may stop at any state $z$, and we assign to its endpoint the \emph{terminal charge} $\dtp{w,z}+\beta_w$, where $\beta_w>0$ is fixed. The \emph{ceiling} of a path is the largest among its terminal charge and the unused levels of its steps. The empty path, which stops immediately at $x$, has ceiling $\dtp{w,x}+\beta_w$. A nonempty path has a strictly smaller ceiling if and only if it ends at a state with a smaller weighted count and all its unused levels are below $\dtp{w,x}+\beta_w$. Hence, $\beta_w$ is the amount by which the unused levels may exceed the starting weighted count while the path prepares the decrease.

\begin{lemma}[Translation of a path]\label{lem:translation}
Let $\pi$ be a path from $x$ to $z$ and let $v\in\nat^d$. Adding $v$ to every state of $\pi$ gives a path from $x+v$ to $z+v$ with the same reactions, which we call the \emph{translation} of $\pi$ by $v$. Its ceiling is the ceiling of $\pi$ plus $\dtp{w,v}$.
\end{lemma}
\begin{proof}
A reaction that is active at a state $x$ is active at every state $x'\ge x$, so the translation is a path. Each of its unused levels and its terminal charge are those of $\pi$ plus $\dtp{w,v}$.
\end{proof}

\begin{definition}[Height and state correction]\label{def:height}
For every state $x\in\nat^d$, define the \emph{height} and the \emph{state correction} by
\begin{equation}\label{eq:height-infimum}
 \begin{split}
 H_w(x)&=\min_{\pi:x\to z}\max\bigl(\{\dtp{w,z}+\beta_w\}\cup\{\dtp{w,x_0-y_1},\ldots,\dtp{w,x_{k-1}-y_k}\}\bigr),\\
 s_w(x)&=\dtp{w,x}-H_w(x).
 \end{split}
\end{equation}
The minimum ranges over all paths $\pi$ from $x$, where $(x_0,y_1\to y_1'),\ldots,(x_{k-1},y_k\to y_k')$ are the steps of $\pi$, with $x_0=x$, and $z$ is the endpoint of $\pi$, including the empty path at $x$, whose ceiling is $\dtp{w,x}+\beta_w$.
\end{definition}

\begin{lemma}[The minimum exists]\label{lem:min-attained}
The minimum in \eqref{eq:height-infimum} is attained.
\end{lemma}
\begin{proof}
At or below the ceiling of the empty path there are only finitely many possible unused levels and terminal charges $\dtp{w,z}+\beta_w$: every integer coordinate with positive weight is bounded, and zero-weight coordinates do not affect these values. The possible path ceilings in this range therefore form a nonempty finite set.
\end{proof}

We call a path from $x$ \emph{optimal} if its ceiling is $H_w(x)$. By Lemma~\ref{lem:min-attained}, there is an optimal path from every state.

\begin{examplecont}
Take $w=(1,0)$ and $\beta_w=2$ in \eqref{eq:example-network}. A path from $(3,0)$ with ceiling below $3$ could only alternate between $(3,0)$ and $(4,1)$. Stopping there has ceiling $5$ and $6$, respectively. At ceiling $3$, it can instead create a $B$ using $0\to B$ and then consume an $A$ and a $B$ using $2A+B\to A$. Repeating reaches $(1,0)$, whose terminal charge is $3$. Hence, $H_w(3,0)=3$. Note that the alternating path visits $(4,1)$, whose weighted count is $4$, although its unused levels are below $3$: the ceiling takes into account the unused levels and the terminal charge, and not the weighted counts of the visited states.
\end{examplecont}

\subsection{Relative orders and the role of weak reversibility}
\begin{definition}[Relative order]\label{def:order}
The \emph{relative order} of a step $(x,y\to y')$ is
\[
 \dtp{w,y}-s_w(x)=H_w(x)-\dtp{w,x-y}.
\]
It is positive, zero, or negative according as the unused level lies below, at, or above the current height.
\end{definition}

\begin{definition}[Preparatory steps and order-zero decreases]\label{def:retained}
We call a step $(x,y\to y')$ \emph{preparatory} if its relative order is nonnegative and $H_w(x-y+y')=H_w(x)$: it preserves height and has unused level at most that height. A step $(x,y\to y')$ of relative order zero with $H_w(x-y+y')<H_w(x)$ is an \emph{order-zero decrease}.
\end{definition}

The source complex of an order-zero decrease at $x$ has order $s_w(x)$, so $T^{s_w(x)}$ is the nominal rate of an order-zero decrease at $x$. The relative order compares the nominal rate $T^{\dtp{w,y}}$ of an active reaction with this rate: a reaction of positive relative order is faster than an order-zero decrease, and a reaction of negative relative order is slower. Theorem~\ref{thm:height} shows that the reactions of positive relative order preserve the height. Hence, the height can increase only along the slower reactions.

A point $x$ in a subset $\mathcal M\subseteq\nat^d$ is \emph{coordinatewise minimal} if there is no $m\in\mathcal M$ with $m\le x$ and $m\ne x$. The finiteness statements of this paper rely on the following classical result. We include its proof for completeness.

\begin{lemma}[Dickson's lemma~\cite{Dickson1913}]\label{lem:dickson}
Every subset of $\nat^d$ has finitely many coordinatewise minimal points, and for every point $x$ of the subset there is a minimal point $m$ with $m\le x$.
\end{lemma}
\begin{proof}
Every infinite sequence of natural numbers has an infinite nondecreasing subsequence. Either one value repeats infinitely often, or every value occurs finitely often and we can successively choose larger values. Apply this observation to each coordinate, retaining a subsequence at each of the finitely many stages. Every infinite sequence in $\nat^d$ therefore has an infinite coordinatewise nondecreasing subsequence. An infinite list of distinct minimal points would contain two with $m\le m'$, contradicting the minimality of $m'$. Finally, for any $x\in\mathcal M$, the set $\{m\in\mathcal M:m\le x\}$ is nonempty and finite. A minimal point in this set is also minimal in $\mathcal M$.
\end{proof}

\begin{theorem}[A bounded correction from path ceilings]\label{thm:height}
For every $\beta_w>0$, the following hold.
\begin{enumerate}[label=(\roman*)]
\item For all $x,x'\in\nat^d$ with $x\le x'$, we have $s_w(x)\le s_w(x')$. Moreover, the correction $s_w$ takes finitely many values, all in $[-\beta_w,\max_{y\in\C}\dtp{w,y}]$.
\item For any linkage class $\mathcal Y$ and any $u\in\nat^d$,
\begin{equation}\label{eq:admissible-height}
 \max\{H_w(y+u)-\dtp{w,u},0\}=\max\{\dtp{w,y}-s_w(y+u),0\}\quad\text{is independent of }y\in\mathcal Y.
\end{equation}
\item Every step of positive relative order is preparatory, and every step of relative order zero is preparatory or an order-zero decrease.
\item Every step $(x,y\to y')$ of positive relative order admits a path back from $x-y+y'$ to $x$: every state on this path has height $H_w(x)$, and every step has the same relative order as $(x,y\to y')$.
\item There is $\gamma_w>0$ such that every step $(x,y\to y')$ satisfies the following:
\[
 \dtp{w,y}-s_w(x)<0\quad\Longrightarrow\quad \dtp{w,y}-s_w(x)\le-\gamma_w,
\]
and, if $(x,y\to y')$ is an order-zero decrease, then $H_w(x-y+y')\le H_w(x)-\gamma_w$.
\end{enumerate}

The gap $\gamma_w$ may depend on $\beta_w$.
\end{theorem}

\begin{proof}
\emph{Monotonicity and boundedness of the correction.}
Let $x,v\in\nat^d$. Lemma~\ref{lem:translation}, applied to an optimal path from $x$, gives
\[
 H_w(x+v)\le H_w(x)+\dtp{w,v},
 \qquad s_w(x+v)\ge s_w(x).
\]
The empty path gives $H_w(x)\le \dtp{w,x}+\beta_w$, hence $s_w(x)\ge-\beta_w$. For the upper bound on $s_w(x)$, we show that every path from $x$ has ceiling at least $\dtp{w,x}-\max_{y\in\C}\dtp{w,y}$. The ceiling of a nonempty path is at least the unused level $\dtp{w,x}-\dtp{w,y}$ of its first step $(x,y\to y')$. The ceiling of the empty path is $\dtp{w,x}+\beta_w$, which is larger than $\dtp{w,x}$ because $\beta_w>0$, and the complex orders are nonnegative. Thus, $H_w(x)\ge\dtp{w,x}-\max_{y\in\C}\dtp{w,y}$, that is, $s_w(x)\le\max_{y\in\C}\dtp{w,y}$. The finiteness of the set of values of $s_w$ is proved at the end, since it relies on the gap in (v).

\emph{Comparison of the heights within a linkage class.}
Let $y$ and $y'$ be two complexes in the same linkage class, and let $u\in\nat^d$. By weak reversibility, there are directed paths of complexes from $y$ to $y'$ and from $y'$ to $y$. Adding $u$ to every complex turns them into paths between the states $y+u$ and $y'+u$, all of whose steps have unused level $\dtp{w,u}$. Following the path from $y+u$ to $y'+u$ by an optimal path from $y'+u$ gives
\[
 H_w(y+u)\le\max\{\dtp{w,u},H_w(y'+u)\},
\]
and the same inequality holds with $y$ and $y'$ exchanged. If the two heights differ, the larger is at most $\dtp{w,u}$. Thus, either $H_w(y+u)=H_w(y'+u)$, or both heights are at most $\dtp{w,u}$. In the first case, $H_w(y+u)-\dtp{w,u}=H_w(y'+u)-\dtp{w,u}$. In the second case, both differences are nonpositive. In both cases,
\[
 \max\{H_w(y+u)-\dtp{w,u},0\}=\max\{H_w(y'+u)-\dtp{w,u},0\}.
\]
Finally, $H_w(y+u)-\dtp{w,u}=\dtp{w,y}-s_w(y+u)$ by the definition of $s_w(y+u)$, which gives the second expression in \eqref{eq:admissible-height}.

\emph{Steps of nonnegative relative order.}
Let $(x,y\to y')$ be a step and let $u=x-y$, so that the relative order of the step is $H_w(y+u)-\dtp{w,u}$. The complexes $y$ and $y'$ are in the same linkage class. If the relative order is positive, then $H_w(y+u)>\dtp{w,u}$, so the two heights are not both at most $\dtp{w,u}$. Hence, $H_w(y'+u)=H_w(y+u)$ and the step is preparatory. If the relative order is zero, then $H_w(y+u)=\dtp{w,u}$, and the inequality above with $y$ and $y'$ exchanged gives $H_w(y'+u)\le H_w(y+u)$. Hence, the step is preparatory or an order-zero decrease. This proves (iii).

Assume again that the relative order of $(x,y\to y')$ is positive. By weak reversibility, there is a directed path of complexes from $y'$ to $y$, and adding $x-y$ to every complex turns it into a path from $x-y+y'$ to $x$. Every state of this path is of the form $x-y+\tilde y$, where $\tilde y$ is a complex in the linkage class of $y$. By \eqref{eq:admissible-height} with $u=x-y$, $\max\{H_w(x-y+\tilde y)-\dtp{w,x-y},0\}=H_w(x)-\dtp{w,x-y}>0$, so $H_w(x-y+\tilde y)=H_w(x)$. Every step of the path has unused level $\dtp{w,x-y}$, so its relative order is $H_w(x)-\dtp{w,x-y}$. This proves (iv).

\emph{A gap for unused levels and terminal charges.}
For each number $c$ in the finite list $\{\dtp{w,y}:y\in\C\}\cup\{-\beta_w\}$, consider ${\mathcal M}_c=\{x\in\nat^d:s_w(x)>c\}$. If nonempty, this set has finitely many minimal points by Lemma~\ref{lem:dickson}. Their positive differences $s_w(x)-c$ have a positive minimum, which remains a lower bound throughout ${\mathcal M}_c$ by the monotonicity of $s_w$. Taking the minimum over the nonempty sets gives $\gamma_w>0$ such that
\begin{equation}\label{eq:height-positive-gap}
 s_w(x)>c\quad\Longrightarrow\quad s_w(x)\ge c+\gamma_w
 \qquad(c\in\{\dtp{w,y}:y\in\C\}\cup\{-\beta_w\}).
\end{equation}
If all these sets are empty, take $\gamma_w=1$. For a step $(x,y\to y')$, the identity $\dtp{w,x-y}-H_w(x)=s_w(x)-\dtp{w,y}$ gives the gap for negative relative orders. For an order-zero decrease $(x,y\to y')$, its unused level $\dtp{w,x-y}$ is $H_w(x)>H_w(x-y+y')$, so
\[
 H_w(x)-H_w(x-y+y')=s_w(x-y+y')-\dtp{w,y'}\ge\gamma_w.
\]
This proves (v).

\emph{Finitely many correction values.} Fix an optimal path $\pi_x$ from every state $x$, and let $F_x$ be the finite set of numbers $\dtp{w,x}-e$, where $e$ ranges over the unused levels of $\pi_x$ and its terminal charge. We claim that
\[
 x\le x',\quad s_w(x')<s_w(x)+\gamma_w\quad\Longrightarrow\quad s_w(x')\in F_x.
\]
Let $v=x'-x$. Translating $\pi_x$ by $v$ gives a path from $x'$ of ceiling $H_w(x)+\dtp{w,v}=H_w(x')+s_w(x')-s_w(x)<H_w(x')+\gamma_w$. Read the unused levels of the translated path in order, with the terminal charge last. At least one of these numbers is at least $H_w(x')$, since every path from $x'$ has at least that ceiling. Take the first such entry. The height there is still $H_w(x')$: every preceding step had unused level below $H_w(x')$, hence positive relative order, and preserved the height by (iii). The excess of this entry above the height is $s_w(z)-\dtp{w,y}$ for a step $(z,y\to y')$, or $s_w(z)+\beta_w$ for the terminal charge. A positive excess would be at least $\gamma_w$ by \eqref{eq:height-positive-gap}, with $c=\dtp{w,y}$ or $c=-\beta_w$ respectively, more than the path allows. The entry therefore equals $H_w(x')$. If $e$ was that entry on $\pi_x$, then $H_w(x')=e+\dtp{w,v}$ and $s_w(x')=\dtp{w,x}-e\in F_x$.

For $k\in\mathbb Z$, let
\[
 \mathcal B_k=\{x\in\nat^d:k\gamma_w\le s_w(x)<(k+1)\gamma_w\}.
\]
Since $s_w$ is bounded, as shown in the first part of the proof, $\mathcal B_k$ is nonempty only for finitely many $k$. By Lemma~\ref{lem:dickson}, each $\mathcal B_k$ has finitely many minimal points. Hence, the union $\mathcal F$ of the sets $F_m$, over the minimal points $m$ of the sets $\mathcal B_k$, is finite. Let $x'\in\nat^d$ and let $k$ be such that $x'\in\mathcal B_k$. By Lemma~\ref{lem:dickson}, there is a minimal point $m$ of $\mathcal B_k$ with $m\le x'$. Since $s_w(m)\ge k\gamma_w$ and $s_w(x')<(k+1)\gamma_w$, we have $s_w(x')<s_w(m)+\gamma_w$. By the claim, $s_w(x')\in F_m\subseteq\mathcal F$.

\end{proof}

\begin{remark}[Finitely many relative orders]\label{rem:finite-orders}
Since $\C$ is finite and $s_w$ takes finitely many values by Theorem~\ref{thm:height}(i), the relative orders $\dtp{w,y}-s_w(x)$ take finitely many values.
\end{remark}

\begin{lemma}[Choosing the terminal charge]\label{lem:terminal-choice}
One can choose $\beta_w\ge2$ and a threshold $R_w\ge\beta_w$ so that, for every $x\in\Gamma$,
\begin{equation}\label{eq:terminal-improvement}
 \dtp{w,x}>R_w-\beta_w\quad\Longrightarrow\quad H_w(x)<\dtp{w,x}+\beta_w.
\end{equation}
\end{lemma}
\begin{proof}
Consider the states $x$ of $\Gamma$ for which some state $z\le x$ of $\Gamma$ has lower weighted count,
\[
 \mathcal D_w=\{x\in\Gamma:\text{some }z\in\Gamma\text{ satisfies }z\le x\text{ and }\dtp{w,z}<\dtp{w,x}\}.
\]
By Lemma~\ref{lem:dickson}, $\mathcal D_w$ has finitely many minimal elements $m_1,\ldots,m_N$. For each $i$, the definition of $\mathcal D_w$ gives a state $z_i\in\Gamma$ with $z_i\le m_i$ and $\dtp{w,z_i}<\dtp{w,m_i}$. Choose a path $\pi_i$ from $m_i$ to $z_i$, which exists because $\Gamma$ is a communicating class. Choose $\beta_w\ge2$ such that, for each $i$, every unused level of $\pi_i$ is smaller than $\dtp{w,m_i}+\beta_w$. If $\mathcal D_w$ is empty, choose $\beta_w=2$.

For every $x\in\mathcal D_w$ there is $i$ with $m_i\le x$. Translating $\pi_i$ by $x-m_i$ gives a path from $x$ with every unused level strictly below $\dtp{w,x}+\beta_w$. Its terminal charge is also strictly below this value, because $\dtp{w,z_i}<\dtp{w,m_i}$. Hence, $H_w(x)<\dtp{w,x}+\beta_w$. Set $R_w$ to $\beta_w$ plus the largest weighted count among the minimal states of $\Gamma$, which are finitely many by Lemma~\ref{lem:dickson}. If $\dtp{w,x}>R_w-\beta_w$, a minimal state $m\le x$ of $\Gamma$ satisfies $\dtp{w,m}<\dtp{w,x}$, so $x\in\mathcal D_w$. This proves \eqref{eq:terminal-improvement}. Note that the height itself may still be below $R_w$.
\end{proof}

Fix $\beta_w,R_w$ as in Lemma~\ref{lem:terminal-choice} and take the gap $\gamma_w$ for this $\beta_w$ from Theorem~\ref{thm:height}(v). Every terminal charge is at least $\beta_w$, so $H_w(x)\ge\beta_w\ge2$.

\subsection{A boundary example}\label{subsec:example}
\begin{figure}[htbp]
\centering
\begin{tikzpicture}[x=1.15cm,y=0.65cm,>=stealth]
 \draw[->] (-.15,0)--(4.8,0) node[right] {$k$};
 \draw[->] (0,-.1)--(0,2.8) node[above] {$s_w(n,k)$};
 \draw[densely dotted,gray] (0,1)--(4.5,1);
 \draw[densely dotted,gray] (0,2)--(4.5,2);
 \node[left] at (0,1) {$1$}; \node[left] at (0,2) {$2$};
 \draw[thick] (0,0)--(1,0)--(1,1)--(2,1)--(2,2)--(4.5,2);
 \foreach \b/\s in {0/0,1/1,2/2,3/2,4/2}{\fill (\b,\s) circle (2pt);\node[below] at (\b,0) {$\b$};}
\end{tikzpicture}
\caption{The correction at any fixed $n\ge5$ in \eqref{eq:example-height}. The dots give its values at integer $B$ counts.}
\label{fig:height}
\end{figure}
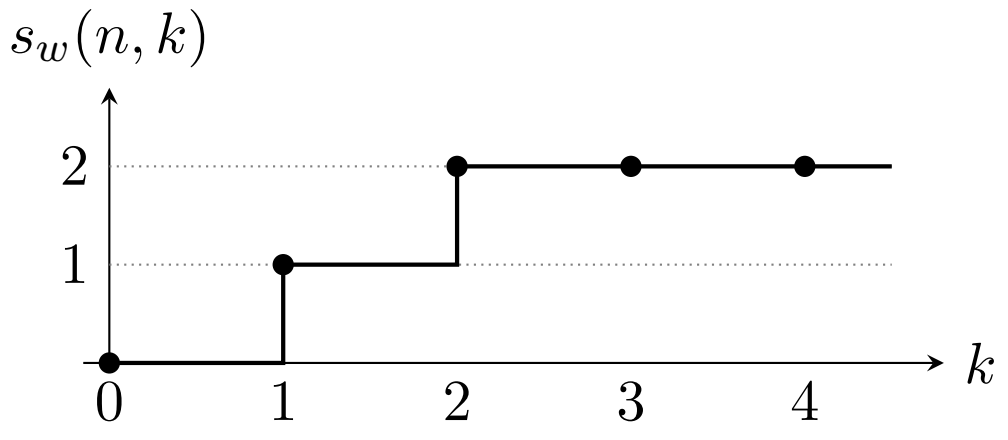

\begin{examplecont}
In \eqref{eq:example-network}, set all four rate constants to one and take $w=(1,0)$. Starting with $n\ge1$, no reaction can remove the last $A$. From any such state, create $B$ when needed to reduce the $A$ count to one, then remove any remaining $B$. This gives a path to $(1,0)$ and, by weak reversibility, a path back. Thus, $n\ge1$ is a single class. Its minimal state is $(1,0)$, and the unique minimal state in $\mathcal D_w$ is $(2,0)$. From $(2,0)$, create a $B$ and use $2A+B\to A$ to reach $(1,0)$; the largest unused level is $2$. Thus, $\beta_w=2$ meets the requirement of Lemma~\ref{lem:terminal-choice} and gives $R_w=3$.

We claim that
\begin{equation}\label{eq:example-height}
 H_w(n,k)=\max\{3,n-2,n-k\},\qquad
 s_w(n,k)=\min\{n-3,2,k\}.
\end{equation}
Every terminal charge is at least $3$, giving the first lower bound. The first step of any nonempty path uses at most two $A$ molecules, leaving at least $n-2$ unused. The empty path has ceiling $n+2$. For the third bound, only $0\to B$ can lower the difference between the current $A$ and $B$ counts. At its first use, at least $n-k$ molecules of $A$ are present and all are unused. If the path never uses it, the endpoint still has at least $n-k$ molecules of $A$.

For the reverse inequality, let $h=\max\{3,n-2,n-k\}$. While the $A$ count exceeds $h$, use $2A+B\to A$. There are enough $B$ molecules because $n-h\le k$, and the unused levels are at most $n-2$. Once the $A$ count is at most $h$, create a $B$ whenever needed and continue the removals until one $A$ remains. Finally remove any remaining $B$, reaching $(1,0)$ with terminal charge $3$. Every step and this terminal charge fit under $h$, proving \eqref{eq:example-height}.

At $n\ge5$, where $n-3\ge2$, the correction is $0$, $1$, or $2$ according as $k=0$, $k=1$ or $k\ge2$; see Figure~\ref{fig:height}. The reactions $2A+B\to A$, $A\to2A+B$, $B\to0$, and $0\to B$ have actual rates $n(n-1)k$, $n$, $k$, and $1$, respectively. For $n\ge6$, where also $H_w(n-1,k-1)=n-3$ in the last row, the table lists the preparatory steps and order-zero decreases. Each destination is specified by its new $B$ count, while its $A$ count changes as prescribed by the reaction.
\[
\begin{array}{c@{\quad}l@{\quad}c@{\quad}l@{\quad}c}
 k&\text{reaction}&\text{new }k&\text{role}&\text{relative order}\\
 \hline
 0&A\to2A+B&1&\text{preparatory}&1\\
 0&0\to B&1&\text{order-zero decrease}&0\\
 1&2A+B\to A&0&\text{preparatory}&1\\
 1&A\to2A+B&2&\text{preparatory}&0\\
 2&2A+B\to A&1&\text{preparatory}&0\\
 k\ge3&2A+B\to A&k-1&\text{order-zero decrease}&0
\end{array}
\]
In the last row, the relative order is $2-s_w(n,k)=0$ and
\[
 H_w(n-1,k-1)=n-3=H_w(n,k)-1,
\]
so the step is an order-zero decrease. A step can also preserve height with negative relative order, as $0\to B$ does from $(n,2)$ to $(n,3)$; such steps will be estimated separately in the drift proof. At $k=2$, the step $2A+B\to A$ has relative order zero, since its complex order $2$ equals $s_w(n,2)$. Its actual rate is $2n(n-1)$.

Finally, start at $(n,1)$ with $n\ge5$. Here $H_w(n,1)=n-1$. Both reactions involving $A$ preserve $H_w$. The reactions $B\to0$ and $0\to B$ change $H_w$ by $+1$ and $-1$, respectively, each at rate one. Consequently
\begin{equation}\label{eq:example-unadjusted}
 \frac{\L e^{H_w}(n,1)}{e^{H_w(n,1)}}=e+e^{-1}-2>0.
\end{equation}
Hence, $e^{H_w}$ is not a Lyapunov function. The height identifies a path to a decrease, but it does not change during the preparatory steps, and the drift estimate must account for the progress made by these steps.
\end{examplecont}

\Needspace{11\baselineskip}
\section{Measuring the preparation before a decrease}\label{sec:races}

In the running example, at $(n,1)$ there is a path that lowers the height, and yet $e^{H_w}$ has positive drift by \eqref{eq:example-unadjusted}. In this section we correct the height with a bounded function that measures the preparation still needed before a decrease, and we use weak reversibility to bound the amount by which a competing reaction can increase it. We keep the class, the direction, and the choices of Section~\ref{sec:height} fixed. All states in this section are in $\Gamma$. Since $\Gamma$ is closed, every path from a state of $\Gamma$ stays in $\Gamma$.

\subsection{Preparations and their translations}
\begin{definition}[Preparation]\label{def:preparation}
A state $x$ is \emph{above the threshold} if $H_w(x)>R_w$. Such a state is \emph{marked} if an order-zero decrease is active there. A \emph{preparation} from a state $x$ is a path of preparatory steps from $x$ to a marked state. The empty path is a preparation whenever $x$ is marked. The final decrease is not part of the preparation.
\end{definition}

Preparations need not be unique, and they may start at or pass through marked states.

\begin{lemma}[Reaching a first decrease]\label{lem:first-decrease}
For every state $x$ with $H_w(x)>R_w$, there is a preparation from $x$.
\end{lemma}
\begin{proof}
Take an optimal path from $x$ to $z$, with steps $(x_0,y_1\to y_1'),\ldots,(x_{k-1},y_k\to y_k')$, where $x_0=x$, and let $x_k=z$. We show that $H_w(z)<H_w(x)$, so that the height changes along the path, and then that the steps before the first change form a preparation. The empty path at $z$ gives $H_w(z)\le\dtp{w,z}+\beta_w$. Moreover, $\dtp{w,z}+\beta_w$ is the terminal charge of the optimal path, so it is at most its ceiling $H_w(x)$. Hence, $H_w(z)\le \dtp{w,z}+\beta_w\le H_w(x)$. If $H_w(z)=H_w(x)$, both inequalities are equalities, so $\dtp{w,z}=H_w(x)-\beta_w>R_w-\beta_w$. Since $z\in\Gamma$, \eqref{eq:terminal-improvement} gives the contradiction $H_w(z)<\dtp{w,z}+\beta_w=H_w(x)$. Thus, $H_w(z)<H_w(x)$.

Let
\[
 j=\min\{i:H_w(x_i)\ne H_w(x)\},
\]
which exists because $H_w(x_k)<H_w(x)$, and satisfies $j\ge1$. For every $i\le j$ we have $H_w(x_{i-1})=H_w(x)$, and the unused level $\dtp{w,x_{i-1}-y_i}$ is at most the ceiling $H_w(x)$ of the path. Hence, the relative order of the step $(x_{i-1},y_i\to y_i')$ satisfies
\[
 H_w(x_{i-1})-\dtp{w,x_{i-1}-y_i}=H_w(x)-\dtp{w,x_{i-1}-y_i}\ge0\qquad(1\le i\le j).
\]
If $i<j$, then $H_w(x_i)=H_w(x_{i-1})$, so the step is preparatory. For $i=j$ we have $H_w(x_j)\ne H_w(x_{j-1})$, so the step is not preparatory and, by Theorem~\ref{thm:height}(iii), its relative order is not positive. By the display above it is nonnegative, hence zero, and Theorem~\ref{thm:height}(iii) shows that the step is an order-zero decrease. Since $H_w(x_{j-1})=H_w(x)>R_w$, the state $x_{j-1}$ is marked, and the first $j-1$ steps form a preparation from $x$.
\end{proof}

\begin{examplecont}
At $w=(1,0)$ and $x=(n,1)$ with $n\ge5$, a preparation and its final decrease are
\[
 (n,1)\xrightarrow{2A+B\to A}(n-1,0)
 \xrightarrow{0\to B}(n-1,1).
\]
The unused levels are $n-2$ and $n-1=H_w(x)$: the first reaction step is the preparation, of relative order one, and the second is the order-zero decrease.
\end{examplecont}

The next lemma reuses a preparation from a state $m$ at a larger state $x\ge m$, by adding the stock $v=x-m$ to every state of the path. Every unused level then increases by $\dtp{w,v}$. The height of the starting state increases by at most $\dtp{w,v}$, by the monotonicity of $s_w$ in Theorem~\ref{thm:height}(i), and exactly by $\dtp{w,v}$ if $s_w(x)=s_w(m)$. We show in Lemma~\ref{lem:preparation-translation} that in this case the translated path keeps the relative orders of the original one until the first decrease of the height, which may occur earlier than in the original path.

\begin{lemma}[Translating a preparation]\label{lem:preparation-translation}
Suppose $m\le x$, $H_w(m)>R_w$ and $s_w(m)=s_w(x)$. Append an order-zero decrease to any preparation from $m$, translate by $x-m$, and stop immediately before the first height decrease. This gives a preparation from $x$: an initial segment of the translated original preparation, with the same relative orders. Any earlier decrease replaces an originally order-zero preparatory step.
\end{lemma}

\begin{proof}
Let $v=x-m$ and let $z$ be the endpoint after the appended decrease. Since $s_w(x)=s_w(m)$, we have $H_w(x)=\dtp{w,m}+\dtp{w,v}-s_w(m)=H_w(m)+\dtp{w,v}$, so $H_w(x)\ge H_w(m)>R_w$. The translated path starts at height $H_w(x)$ and ends at height
\[
 H_w(z+v)\le H_w(z)+\dtp{w,v}<H_w(m)+\dtp{w,v}=H_w(x).
\]
The first inequality holds because $s_w(z)\le s_w(z+v)$ by Theorem~\ref{thm:height}(i), and the second because the preparation keeps the height $H_w(m)$ and the appended decrease lowers it. Hence, there is a first height change along the translated path. The original steps start at height $H_w(m)$ and have nonnegative relative order, so their unused levels are at most $H_w(m)$, and the translated unused levels are at most $H_w(x)$. Thus, the step of the first height change has nonnegative relative order and is not preparatory, so by Theorem~\ref{thm:height}(iii) it is an order-zero decrease. Its preceding steps form a preparation. Along those steps, both the height and the unused level increase by $\dtp{w,v}$ under translation, preserving their relative orders. The first decrease has unused level $H_w(x)$, so its original relative order is zero.
\end{proof}

\subsection{A bounded correction for preparation}
In Section~\ref{sec:height} the height was defined by minimizing the ceiling of a path, which is a maximum over its steps. We now keep the height fixed and assign to a preparation a cost that sums contributions of its steps and of the decrease that follows it (Definition~\ref{def:cost}). By Remark~\ref{rem:finite-orders}, the relative orders take finitely many values. Let $\mathcal P_w$ be the set of the nonnegative relative orders of the steps starting above the threshold.

\begin{definition}[Cost of preparing a decrease]\label{def:cost}
Let $x$ be above the threshold. The steps of a preparation from $x$ start at states of height $H_w(x)$. The \emph{cost} of such a step of relative order $p$ is $H_w(x)^{-p}$, and the \emph{cost} of the preparation is one plus the sum of the costs of its steps. The one accounts for the order-zero decrease that follows the preparation, whose relative order is zero. A larger relative order gives a smaller cost. Let $\pathcost_w(x)$ be the infimum of the costs of the preparations from $x$.
\end{definition}

\begin{lemma}[A finite menu of path costs]\label{lem:route-distance}
There are nonnegative integers $\ell_w,K_w$ such that, at every $x$ above the threshold, some preparation of length at most $\ell_w$ attains $\pathcost_w(x)$. Moreover, for every $x$ above the threshold, the set $\{\pathcost_w(x'):H_w(x')=H_w(x)\}$ has at most $K_w+1$ elements.
\end{lemma}

\begin{proof}
If no state lies above the threshold, take $K_w=\ell_w=0$. Otherwise, for each of the finitely many values $s_*$ taken by the correction above the threshold (Theorem~\ref{thm:height}(i)), consider the pairs $(x,\mu)\in\Gamma\times\nat^{|\mathcal P_w|}$ such that $H_w(x)>R_w$, $s_w(x)=s_*$, and some preparation from $x$ has count vector $\mu$, where $\mu_p$ counts steps of relative order $p$. The cost of such a preparation is $1+\sum_p \mu_p H_w(x)^{-p}$. Hence, if two preparations from the same state have count vectors $\mu'\le\mu$, then the cost of the first is not larger than the cost of the second. This is the reason why we consider the pair $(x,\mu)$ and not only $x$. Choose a preparation at each of the finitely many minimal pairs given by Lemma~\ref{lem:dickson}.

Let $x$ be above the threshold. By Lemma~\ref{lem:first-decrease}, there is a preparation from $x$. Given any preparation from $x$ with counts $\mu$, choose a minimal pair $(m,\mu^0)\le(x,\mu)$ with the same correction. Lemma~\ref{lem:preparation-translation} translates its chosen preparation into one from $x$ with counts $\mu'\le \mu^0\le \mu$. At the height of $x$, these counts cost no more than $\mu$. The translated preparations have count vectors drawn from the prefixes of finitely many chosen paths. At $x$, keep just the vectors on this list realized by such translations. Every preparation from $x$ is matched by a choice on this list of no greater cost, and each listed choice is itself realized by a preparation from $x$. Minimizing over the list therefore gives exactly $\pathcost_w(x)$. The longest chosen preparation bounds the length needed. At any fixed height each count vector gives only one cost, so we may take for $K_w$ the number of count vectors on the list.
\end{proof}

\begin{examplecont}
Let $(n,k)$ be above the threshold. The minimum cost $\pathcost_w(n,k)$ is $1$ for $k=0$ or $k\ge3$, $1+1/H_w(n,k)$ for $k=1$, and $2+1/H_w(n,k)$ for $k=2$. The one accounts for the decrease that follows the preparation, and a preparatory step of relative order one costs $1/H_w(n,k)$.

At $(n,1)$ we have $H_w(n,1)=n-1$. The two preparatory reactions have rates $n(n-1)$ and $n$ and change the minimum cost by $-1/(n-1)$ and $+1$. Their combined drift is
\[
 -n(n-1)/(n-1)+n=0.
\]
Hence, the minimum cost has zero drift along these two reactions, although the first one is much faster: the faster step lowers the cost by $1/(n-1)$ only, while the slower step raises it by one. Now multiply the first gap between consecutive minimum costs, from $1$ to $1+1/(n-1)$, by one half, and the next gap, from $1+1/(n-1)$ to $2+1/(n-1)$, by one quarter. The two changes become $-1/(2(n-1))$ and $+1/4$, and with the original rates the combined drift becomes $-n/2+n/4<0$. The construction below discounts successive gaps in this way at every height, with a margin that absorbs the error in our estimates of the rates.
\end{examplecont}

\subsection{Progress and detours at the actual rates}\label{sec:local}
A preparatory reaction may lead away from a cheapest preparation. Weak reversibility bounds how much this detour can add to the remaining cost.

\begin{lemma}[Cost of a detour]\label{lem:detour}
For each reaction choose a return path in the reaction graph, and let $\Lambda\ge1$ bound their lengths. If $(x,y\to y')$ is a preparatory step of relative order $p$ and $H_w(x)>R_w$, then
\[
 \pathcost_w(x-y+y')\le \pathcost_w(x)+\Lambda H_w(x)^{-p}.
\]
\end{lemma}
\begin{proof}
Carry the unused stock along the chosen return path. If $p>0$, every return step preserves height and has order $p$ by Theorem~\ref{thm:height}(iv). Returning to $x$ and following a cheapest preparation from $x$ gives the bound. If $p=0$, the unused level stays at $H_w(x)$. By Theorem~\ref{thm:height}(iii), either we return to $x$ through states of height $H_w(x)$, along at most $\Lambda$ preparatory steps of cost one, or some return step is a first decrease. In the second case the return steps before it form a preparation from $x-y+y'$ of at most $\Lambda$ steps, so $\pathcost_w(x-y+y')\le1+\Lambda\le\pathcost_w(x)+\Lambda$, because every cost is at least one. This is a bound on the cheapest continuation from $x-y+y'$; the chain need not follow the return path.
\end{proof}

In Example~\ref{ex:running} with unit rate constants and $n\ge5$, $(n,1)\to(n-1,0)$ and its reverse have rates $n(n-1)$ and $n-1$, yet both steps have relative order one and cost $1/(n-1)$. To estimate the drift at a given state, we next compare, for the reactions that are active there, the product of the rate and the cost of the step.

Recall that $T=\sum_i(x_i+1)$. Let $M\ge1$ be the largest number of molecules in a source complex. As in the example, the gaps between consecutive minimum costs at the height $H_w(x)$ are discounted by successive powers of $H_w(x)^{-\sigma_w}$, with
\[
 \sigma_w=\frac{\gamma_w}{K_w+1}.
\]
There are at most $K_w$ gaps above the first cost, so the largest discount is $H_w(x)^{-K_w\sigma_w}$, and the exponent $K_w\sigma_w=\gamma_w-\sigma_w$ leaves a margin $\sigma_w$ below $\gamma_w$, which absorbs the errors in the estimates of the rates.

\begin{definition}[A fitting direction]\label{def:fit}
Choose $0<E_w<1/2$ so that $2ME_w<\sigma_w$. We say that $w$ \emph{fits} at $x$ if $\|\log_T (x+1)-w\|_\infty<E_w$.
\end{definition}

Fitting compares the counts with the total size $T$, whereas the cost of a preparatory step is expressed in terms of the height $H_w(x)$. At a large fitting state these two scales are comparable. Indeed, a largest coordinate $(x_i+1)\ge T/d$ has $\log_T (x_i+1)\to1$ and hence $w_i\ge1/3$ eventually. Since $s_w$ is bounded by Theorem~\ref{thm:height}(i), we have $H_w(x)=\dtp{w,x+1}+O(1)$. Hence, $H_w(x)$ lies between a positive multiple of $T$ and $T+O(1)$. In particular $H_w(x)>R_w$.

\begin{lemma}[Comparison with the actual rates]\label{lem:actual-rates}
There is a constant $C$ such that, if $w$ fits at a sufficiently large state $x$, each reaction $y\to y'$ active at $x$ satisfies
\begin{equation}\label{eq:actual-relative-rate}
 C^{-1}H_w(x)^{s_w(x)-ME_w}\le\lambda_{y\to y'}(x)H_w(x)^{s_w(x)-\dtp{w,y}}\le CH_w(x)^{s_w(x)+ME_w}.
\end{equation}
\end{lemma}
\begin{proof}
Every factor in $x_i(x_i-1)\cdots(x_i-y_i+1)$ lies between $(x_i+1)/(y_i+1)$ and $(x_i+1)$: writing $x_i=y_i+t$, even the smallest factor $t+1$ is at least $(x_i+1)/(y_i+1)$. The case $y_i=0$ contributes the empty product. Thus, $\lambda_{y\to y'}(x)$ is between fixed positive multiples of $\prod_i (x_i+1)^{y_i}$, since $\kappa_{\min}\le\kappa_{y\to y'}(x)\le\kappa_{\max}$. Fitting bounds that product between $T^{\dtp{w,y}-ME_w}$ and $T^{\dtp{w,y}+ME_w}$. We may replace $T$ by $H_w(x)$ up to constant factors. In particular, a species absent at $x$ causes no difficulty for an active reaction: its source coefficient is zero.
\end{proof}

\begin{definition}[Preparation correction]\label{def:preparation-correction}
Let $\delta_w=\gamma_w/[2(1+\ell_w)]$. Set $D_w(x)=0$ if $H_w(x)\le R_w$. Let $x$ be above the threshold, and let $\tau_0<\cdots<\tau_m$ be the distinct values of $\pathcost_w(x')$ over the states $x'$ with $H_w(x')=H_w(x)$, where $m\le K_w$ by Lemma~\ref{lem:route-distance}. Let $\tau_{-1}=0$, let $k$ be such that $\pathcost_w(x)=\tau_k$, and discount successive gaps:
\begin{equation}\label{eq:preparation-correction}
 D_w(x)=\delta_w\sum_{i=0}^{k}H_w(x)^{-i\sigma_w}(\tau_i-\tau_{i-1}).
\end{equation}
\end{definition}
When the minimum cost changes from $\tau_k$ to $\tau_j$, the correction $D_w$ changes by the sum of the discounted gaps between these two costs. If the cost decreases, these gaps have index at most $k$, and if it increases, they have index at least $k+1$. Hence, every gap involved in an increase carries at least one more factor $H_w(x)^{-\sigma_w}$ than the gap of index $k$. Above the threshold, $D_w(x)\le\delta_w\pathcost_w(x)\le\delta_w(1+\ell_w)=\gamma_w/2$, since, by Lemma~\ref{lem:route-distance}, the minimum cost is attained with at most $\ell_w$ preparatory steps, each of which costs at most one.

When we combine several directions in Section~\ref{sec:selection}, a penalty $b(H_w(x))$, which depends on the direction, is subtracted from the height $H_w(x)$. On a first reading, one may take $b=0$.

\begin{definition}[Penalty and score]\label{def:score}
A \emph{penalty} is a function $b$ such that an increase of the height increases the penalty by a nonnegative amount, which is at most half of the increase of the height:
\begin{equation}\label{eq:half-charge}
 0\le b(v)-b(\tilde u)\le (v-\tilde u)/2\qquad(2\le \tilde u\le v).
\end{equation}
We write $G_w(x)=H_w(x)-b(H_w(x))$ for the score.
\end{definition}

A preparatory step leaves the score fixed, and an order-zero decrease lowers it by at least $\gamma_w/2$, by Theorem~\ref{thm:height}(v) and \eqref{eq:half-charge}.

\begin{lemma}[Full drift from preparation]\label{lem:uniform-preparation}
For the correction $D_w$ of Definition~\ref{def:preparation-correction} and any penalty satisfying \eqref{eq:half-charge}, set $W_w=e^{G_w}(1+D_w)$. There is a constant $C_w$ such that every step $(x,y\to y')$ satisfies
\begin{equation}\label{eq:one-jump}
 W_w(x-y+y')\le C_wW_w(x).
\end{equation}
Moreover, there is $c_w>0$ such that, at every sufficiently large state $x$ at which $w$ fits,
\begin{equation}\label{eq:preparation-drift}
 \L W_w(x)\le-c_w T^{-\nu_w}W_w(x),\qquad \nu_w=\beta_w+\gamma_w.
\end{equation}
\end{lemma}

\begin{proof}
\emph{Proof of \eqref{eq:one-jump}.} Since $s_w$ is bounded by Theorem~\ref{thm:height}(i), there is a constant $A_w$ such that $|H_w(x-y+y')-H_w(x)|\le A_w$ for every step $(x,y\to y')$. By \eqref{eq:half-charge}, we also have $|G_w(x-y+y')-G_w(x)|\le A_w$. Since $0\le D_w\le\gamma_w/2$, \eqref{eq:one-jump} holds with $C_w=e^{A_w}(1+\gamma_w/2)$.

\emph{Setting.} Let $x$ be a sufficiently large state at which $w$ fits. Then, $H_w(x)>R_w$ and $H_w(x)$ is comparable to $T$, as shown before Lemma~\ref{lem:actual-rates}. Let $\tau_0<\cdots<\tau_m$ be the distinct minimum costs at the height $H_w(x)$, as in Definition~\ref{def:preparation-correction}, and let $k$ be such that $\pathcost_w(x)=\tau_k$. For every reaction $y\to y'$ active at $x$, let
\[
 \Delta_{y\to y'}=\frac{W_w(x-y+y')-W_w(x)}{e^{G_w(x)}},\qquad\text{so that}\qquad\frac{\L W_w(x)}{e^{G_w(x)}}=\sum_{y\to y'}\lambda_{y\to y'}(x)\Delta_{y\to y'},
\]
where the sum ranges over the reactions active at $x$. Recall that the relative order of the step $(x,y\to y')$ is $p=\dtp{w,y}-s_w(x)$. Hence, \eqref{eq:actual-relative-rate} reads
\begin{equation}\label{eq:rate-times-cost}
 C^{-1}H_w(x)^{s_w(x)-ME_w}\le\lambda_{y\to y'}(x)H_w(x)^{-p}\le CH_w(x)^{s_w(x)+ME_w}.
\end{equation}
By Theorem~\ref{thm:height}(iii) and (v), every step $(x,y\to y')$ is preparatory, or it is an order-zero decrease, or its relative order is at most $-\gamma_w$. We consider the three cases separately.

\emph{Preparatory steps.} Let $(x,y\to y')$ be a preparatory step of relative order $p$. Its destination has height $H_w(x)$, so $\pathcost_w(x-y+y')=\tau_j$ for some $j$, and the score does not change. Hence, $\Delta_{y\to y'}=D_w(x-y+y')-D_w(x)$, and \eqref{eq:preparation-correction} gives
\begin{equation}\label{eq:gap-bounds}
 \Delta_{y\to y'}\le
 \begin{cases}
 -\delta_wH_w(x)^{-k\sigma_w}(\tau_k-\tau_j)&\text{if }j\le k,\\
 \delta_wH_w(x)^{-(k+1)\sigma_w}(\tau_j-\tau_k)&\text{if }j>k.
 \end{cases}
\end{equation}
If $j>k$, then $\tau_j-\tau_k\le\Lambda H_w(x)^{-p}$ by Lemma~\ref{lem:detour}, and \eqref{eq:rate-times-cost} gives
\[
 \lambda_{y\to y'}(x)\Delta_{y\to y'}\le\delta_w\Lambda CH_w(x)^{s_w(x)+ME_w-(k+1)\sigma_w}.
\]

\emph{Order-zero decreases.} Let $(x,y\to y')$ be an order-zero decrease. It lowers the score by at least $\gamma_w/2$, so $W_w(x-y+y')\le e^{G_w(x)-\gamma_w/2}(1+\gamma_w/2)$. Since $W_w(x)\ge e^{G_w(x)}$, we have
\[
 \Delta_{y\to y'}\le-a_w,\qquad\text{where }a_w=1-e^{-\gamma_w/2}(1+\gamma_w/2)>0.
\]

\emph{Steps of negative relative order.} Let $(x,y\to y')$ have relative order $p\le-\gamma_w$. By \eqref{eq:one-jump}, $\Delta_{y\to y'}\le C_w(1+\gamma_w/2)$. Since $\gamma_w=(K_w+1)\sigma_w$ and $k\le K_w$, \eqref{eq:rate-times-cost} gives
\[
 \begin{split}
 \lambda_{y\to y'}(x)\Delta_{y\to y'}&\le C_w(1+\gamma_w/2)CH_w(x)^{s_w(x)+ME_w-\gamma_w}\\
 &\le C_w(1+\gamma_w/2)CH_w(x)^{s_w(x)+ME_w-(k+1)\sigma_w}.
 \end{split}
\]

\emph{A step with a guaranteed decrease.} By Lemma~\ref{lem:route-distance}, $\pathcost_w(x)$ is the cost of a preparation from $x$. If this preparation is empty, let $(x,y_*\to y_*')$ be an order-zero decrease active at $x$; otherwise, let it be the first step of the preparation. Let $p_*$ be its relative order. In the first case, $p_*=0$ and $\Delta_{y_*\to y_*'}\le-a_w$. In the second case, the step is preparatory, and the remaining steps form a preparation from $x-y_*+y_*'$. Hence, $\pathcost_w(x-y_*+y_*')\le\tau_k-H_w(x)^{-p_*}$, and the first case of \eqref{eq:gap-bounds} gives $\Delta_{y_*\to y_*'}\le-\delta_wH_w(x)^{-k\sigma_w}H_w(x)^{-p_*}$. In both cases, \eqref{eq:rate-times-cost} gives
\[
 \lambda_{y_*\to y_*'}(x)\Delta_{y_*\to y_*'}\le-\min\{a_w,\delta_w\}C^{-1}H_w(x)^{s_w(x)-ME_w-k\sigma_w}.
\]

\emph{Conclusion.} By the three cases above, every active reaction other than $y_*\to y_*'$ contributes to the sum at most $C'H_w(x)^{s_w(x)+ME_w-(k+1)\sigma_w}$, where $C'=C\max\{\delta_w\Lambda,C_w(1+\gamma_w/2)\}$. Hence,
\[
 \frac{\L W_w(x)}{e^{G_w(x)}}\le-H_w(x)^{s_w(x)-ME_w-k\sigma_w}\Bigl(\min\{a_w,\delta_w\}C^{-1}-|\Re|C'H_w(x)^{2ME_w-\sigma_w}\Bigr).
\]
Since $2ME_w<\sigma_w$, the quantity in parentheses is at least $\min\{a_w,\delta_w\}/(2C)$ if $H_w(x)$ is large enough. Moreover, $s_w(x)\ge-\beta_w$ by Theorem~\ref{thm:height}(i) and $ME_w+k\sigma_w<(K_w+1)\sigma_w=\gamma_w$, so $H_w(x)^{s_w(x)-ME_w-k\sigma_w}\ge H_w(x)^{-\nu_w}$. Finally, $e^{G_w(x)}\ge W_w(x)/(1+\gamma_w/2)$, and $H_w(x)\le2T$ if $T$ is large enough. This gives \eqref{eq:preparation-drift}.
\end{proof}

\begin{examplecont}
Keep the four rate constants equal to one and take $b=0$. The preparatory component containing $(n,1)$ has exactly three states, all at height $n-1$ for large $n$:
\[
 z_0=(n-1,0),\qquad z_1=(n,1),\qquad z_2=(n+1,2).
\]
Their state corrections are $0,1,2$. The actual rates of the preparatory steps are
\[
 z_0\xrightarrow{n-1}z_1,
 \qquad z_1\xrightarrow{n(n-1)}z_0,\quad z_1\xrightarrow{n}z_2,
 \qquad z_2\xrightarrow{2n(n+1)}z_1.
\]
The order-zero decrease at $z_0$ has rate one. Using the three minimum costs computed above and writing $\vartheta=(n-1)^{-\sigma_w}$, the construction gives
\[
 D_w(z_0)=\delta_w,\qquad D_w(z_1)=\delta_w\Bigl(1+\frac{\vartheta}{n-1}\Bigr),\qquad D_w(z_2)=\delta_w\Bigl(1+\frac{\vartheta}{n-1}+\vartheta^2\Bigr).
\]
Counting the decrease as a jump to a state where $D_w$ vanishes, the drifts of $D_w$ at $z_0,z_1,z_2$, divided by $\delta_w$, are
\[
 -1+\vartheta,\qquad -n\vartheta(1-\vartheta),\qquad -2n(n+1)\vartheta^2.
\]
All three are negative. The first cost gap pays for the final decrease at $z_0$, despite its fast competitor. The next gaps account for progress during preparation.

Return to $z_1=(n,1)$, where \eqref{eq:example-unadjusted} gave $\L e^{H_w}(n,1)/e^{H_w(n,1)}=e+e^{-1}-2>0$. Here we can take $\gamma_w=1$ and $K_w=2$, so $\vartheta=(n-1)^{-1/3}$. In $\L W_w(n,1)/e^{H_w(n,1)}$, the step $2A+B\to A$ contributes $-\delta_w n\vartheta$, while $A\to2A+B$ contributes $\delta_w n\vartheta^2$. Each of the remaining reactions has rate one and changes the height by one, so boundedness of $D_w$ bounds both contributions. Altogether,
\[
 \frac{\L W_w(n,1)}{e^{H_w(n,1)}}=-\delta_w n\vartheta+\delta_w n\vartheta^2+O(1)<0
\]
for large $n$. The loss is of size $n^{2/3}$ and the competing gain is of size $n^{1/3}$. The correction has changed the sign of the drift at the very state where the height alone failed.
\end{examplecont}

\section{Selecting directions that fit the state}\label{sec:selection}

In this section we choose the finite family $\mathcal A$ of directions for the fixed class, together with their penalties. The goal is Proposition~\ref{prop:selection}: at a large state, every direction whose score is close to the largest one fits. Recall that the score of a direction $w$ at $x$ is
\[
 G_w(x)=H_w(x)-b_w(H_w(x)).
\]
Note that if the score of $w$ is $m\log T$ below the score of another direction, then $e^{G_w}$ is smaller by a factor $T^{-m}$. For this reason, we measure the differences between scores in multiples of $\log T$.

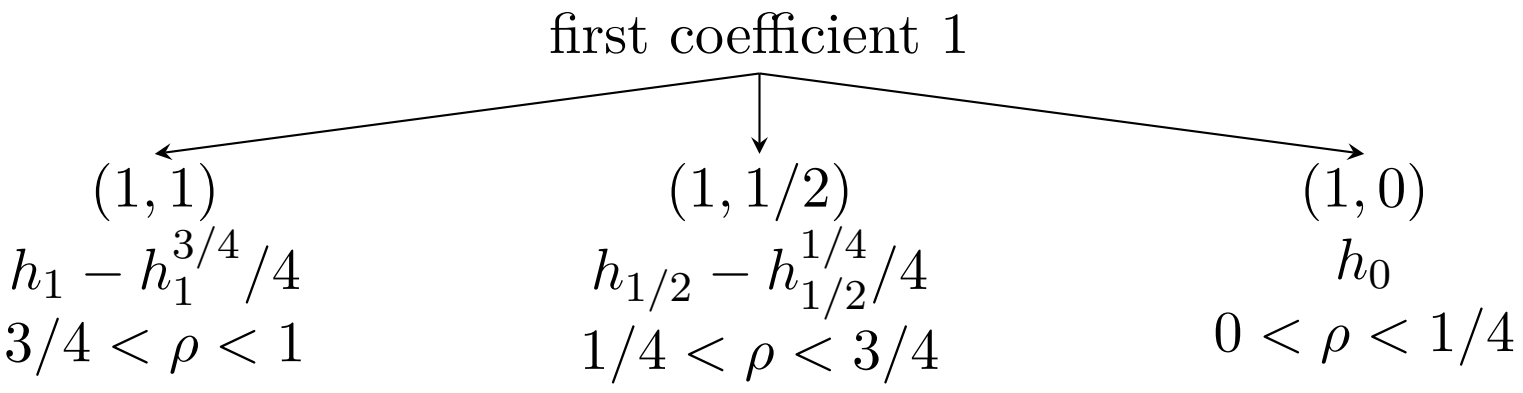
\begin{figure}[htbp]
\centering
\begin{tikzpicture}[x=1cm,y=1cm,>=stealth]
 \node (root) at (0,0) {first coefficient $1$};
 \node[align=center,anchor=north,inner ysep=1.5pt] (first) at (-4,-0.80) {$(1,1)$\\$h_1-h_1^{3/4}/4$\\$3/4<\rho<1$};
 \node[align=center,anchor=north,inner ysep=1.5pt] (middle) at (0,-0.80) {$(1,1/2)$\\$h_{1/2}-h_{1/2}^{1/4}/4$\\$1/4<\rho<3/4$};
 \node[align=center,anchor=north,inner ysep=1.5pt] (last) at (4,-0.80) {$(1,0)$\\$h_0$\\$0<\rho<1/4$};
 \draw[->] (root.south) -- (first.north);
 \draw[->] (root.south) -- (middle.north);
 \draw[->] (root.south) -- (last.north);
\end{tikzpicture}
\caption{At $x+1=(n,\lfloor n^\rho\rfloor)$ the indicated score is the largest of the three for each fixed $\rho$ in its displayed interval and all sufficiently large $n$.}
\label{fig:slope-scores}
\end{figure}

\begin{example}[Choosing the weight of a less abundant species]\label{ex:scores}
Suppose our two-species class contains shifted states $x+1=(n,\lfloor n^\rho\rfloor)$ with $0<\rho<1$ as $n$ grows. The height for $(1,t)$ is $h_t=n+tn^\rho+O(1)$. Raising $t$ from zero to $1/2$ increases the height by $n^\rho/2+O(1)$. If the direction $(1,1/2)$ pays the penalty $h_{1/2}^{1/4}/4$, then its score is larger than the score of $(1,0)$ for $\rho>1/4$ and smaller for $\rho<1/4$, since $h_{1/2}/n\to1$. Figure~\ref{fig:slope-scores} adds a third choice.
\end{example}

\begin{lemma}[Tagged partition]\label{lem:partition}
Let $\bar t>0$ and $\omega:[0,\bar t]\to(0,\infty)$. Then, there are points $0=\tilde u_0<\tilde u_1<\cdots<\tilde u_k=\bar t$ and tags $0=t_1<t_2<\cdots<t_k=\bar t$ such that, for every $1\le j\le k$,
\[
 \tilde u_{j-1}\le t_j\le \tilde u_j,\qquad \tilde u_j-t_j<\omega(t_j),\qquad t_j-\tilde u_{j-1}<\omega(t_j),
\]
and $t_j<\tilde u_j$ for $j<k$.
\end{lemma}
\begin{proof}
Let $\mathcal M$ be the set of $\tilde u\in(0,\bar t)$ for which the statement holds with $\tilde u$ in place of $\bar t$, except that the last tag is only required to satisfy $t_k<\tilde u$. Every $\tilde u<\min\{\bar t,\omega(0)\}$ lies in $\mathcal M$, with the single interval $[0,\tilde u]$ tagged $0$. Let $\tilde u^*=\sup\mathcal M$ and choose $\tilde u\in\mathcal M$ with $\tilde u^*-\tilde u<\omega(\tilde u^*)$. If $\tilde u^*<\bar t$, choose $U\in(\tilde u^*,\bar t)$ with $U-\tilde u^*<\omega(\tilde u^*)$. Adding the interval $[\tilde u,U]$ tagged $\tilde u^*$ shows $U\in\mathcal M$, contradicting the choice of $\tilde u^*$. Hence, $\tilde u^*=\bar t$, and adding $[\tilde u,\bar t]$ tagged $\bar t$ gives the partition of $[0,\bar t]$. Each tag other than the last is smaller than the right endpoint of its interval, which is the left endpoint of the next. The tags therefore increase strictly.
\end{proof}

The selection in Proposition~\ref{prop:selection} adapts the affine approximation of~\cite[Lemma~2.1 and Corollary~2.2]{OpenAIBP2026}, and we give a self-contained proof.

\begin{proposition}[Selecting directions that fit]\label{prop:selection}
Fix the class $\Gamma$ and the accuracies $E_w$ of Definition~\ref{def:fit}. There is a finite family $\mathcal A$ of directions containing $\mathbf1$, with nonnegative penalties $b_w$ satisfying \eqref{eq:half-charge} and $b_w(h)=o(h)$ as $h\to\infty$. For every fixed $m>0$ and all sufficiently large $T$, every $w\in\mathcal A$ satisfying
\[
 G_w(x)\ge\max_{w'\in\mathcal A}G_{w'}(x)-m\log T
\]
fits at $x$ in the sense of Definition~\ref{def:fit}.
\end{proposition}
\begin{proof}
\textbf{Choosing the coefficients.} We construct the directions of $\mathcal A$ one coefficient at a time, in nonincreasing order. The coefficients cannot be chosen independently of each other. For example, for $3A\rightleftharpoons Z+3B$ give $Z,A,B$ the weights $1,p,q$. Both reactions preserve the weighted count when $3p=1+3q$, whatever their rates are. We therefore allow the choice of each coefficient to depend on the preceding ones. The first coefficient is $1$. After a coefficient $\bar t>0$, we partition $[0,\bar t]$ into closed intervals $[L,U]$ that meet only at successive endpoints, and to each interval we assign a \emph{tag} $t\in[L,U]$. The tags are the possible values of the next coefficient, and the interval of a tag collects the exponents that the tag is meant to approximate. The tags $0$ and $\bar t$ always occur, so that a direction can end with zeros or repeat $\bar t$.

A \emph{completion} of a list of coefficients is a direction with $d$ coefficients that starts with that list. The coefficients are appended from left to right, but the accuracy $E_w$ of Definition~\ref{def:fit} is only defined for a full direction $w$. For this reason, the interval of a tag $t$ is chosen \emph{after} the completions of the list ending with $t$ have been constructed. We proceed by induction on the number of positions still to be filled. A list of $d$ coefficients is its own completion. After $t=0$, the remaining positions are filled with zeros. After any $t>0$, the inductive hypothesis gives finitely many completions. We take the smallest accuracy among the completions of $t$ and their permutations, and we require that the interval tagged $t$ lies strictly within half this accuracy of $t$ and avoids each endpoint of $[0,\bar t]$ different from $t$.

Such a partition exists by Lemma~\ref{lem:partition}, with $\omega(t)$ the half accuracy allowed at the tag $t$. By the lemma, the tags increase strictly from $0$ to $\bar t$, every interval is nondegenerate, $L>0$ whenever $t>0$, and $U<\bar t$ whenever $t<\bar t$. The family $\mathcal A$ consists of all the completions obtained in this way, together with all their permutations. Repeating $1$ gives the all-ones direction.

\paragraph{Assigning the penalties.} The initial coefficient $1$ and the zero coefficients carry no penalty. Each later positive coefficient is charged $h^c/(2d)$ at height $h$, where $c$ is the smallest left endpoint among the intervals chosen from the second coefficient up to and including this one. Each such endpoint lies in $(0,1)$: it is positive because the coefficients up to this one are positive and the interval of a positive tag does not start at $0$, and it is below one because every interval of the partition of $[0,1]$ has $L<1$. Since $c$ is a minimum over all the intervals chosen so far, the exponents do not increase along the list, even when a coefficient is repeated. The penalty $b_w(h)$ is the sum of these charges, and it is the same function for all the permutations of $w$.

\begin{figure}[htbp]
\centering
\begin{tikzpicture}[x=15cm,y=1cm,>=stealth]
 \fill[gray!13] (.1,-.1) rectangle (.3,.12);
 \fill[gray!13] (.48,-.1) rectangle (.5,.12);
 \draw[thick] (0,0)--(.5,0);
 \foreach \s in {0,.1,.3,.48,.5} {\draw (\s,-.11)--(\s,.12);}
 \foreach \s/\lab in {0/0,.2/0.2,.4/0.4,.5/0.5} {\fill (\s,0) circle[radius=1.6pt];\node[above=4pt] at (\s,0){$\lab$};}
 \foreach \s/\lab in {.1/0.1,.3/0.3,.48/0.48} {\node[below=4pt] at (\s,0){$\lab$};}
\end{tikzpicture}
\caption{The four possible coefficients are marked above the line and the three cuts below. If the preceding charge had exponent $0.47$, choosing $0.4$ gives exponent $0.3$ and repeating $0.5$ keeps $0.47$, although its own left endpoint is $0.48$.}\label{fig:interval-chain}
\end{figure}
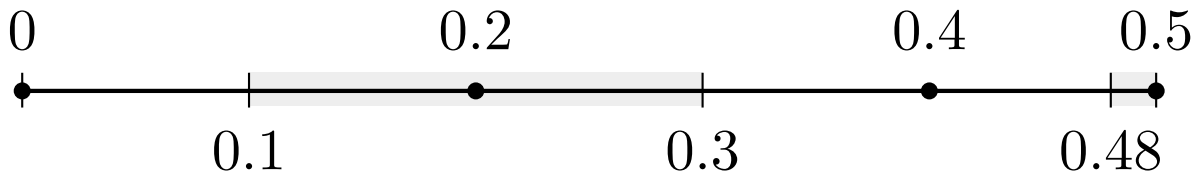

Figure~\ref{fig:interval-chain} illustrates a choice at $\bar t=0.5$: suppose the smallest completion accuracy is $0.08$ at tag $0.5$ and $0.4$ at the other tags. This explains why the last interval is so short. Each charge is $o(h)$ because $c<1$. Since $ch^{c-1}\le1$ for $h\ge1$ and fewer than $d$ terms occur, both $b_w$ and the charge for any initial part satisfy \eqref{eq:half-charge}.

\paragraph{Consequences of a high score.} Since $H_{w'}=\dtp{w',x}-s_{w'}$ with $s_{w'}$ bounded by Theorem~\ref{thm:height}(i), there is $c_0\ge0$ such that $|H_{w'}(x)-\dtp{w',x+1}|\le c_0$ for all $w'\in\mathcal A$ and $x\in\Gamma$. Fix $m>0$, a state $x$ and a direction $w\in\mathcal A$ whose score is within $m\log T$ of the maximum. Let $\varepsilon$ be half the smallest accuracy among the permutations of $w$. Then $2\varepsilon\le E_w$, and the intervals of $w$ lie strictly within $\varepsilon$ of their tags because we included all permutations when choosing their widths. We show that $|\log_T (x_i+1)-w_i|<2\varepsilon$ for every $i$ as soon as $T$ exceeds a bound depending only on $\mathcal A$, $c_0$, and $m$. Below, $C$ stands for constants of this kind.

Order the coordinates by decreasing $w_i$ and, within ties, by decreasing $(x_i+1)$.

\emph{The height is comparable to $T$.} We have $H_w(x)\le \dtp{w,x+1}+c_0\le T+c_0$. The all-ones direction has height $H_{\mathbf1}(x)\ge T-c_0$, and by \eqref{eq:half-charge} its score is at least $H_{\mathbf1}(x)/2-C$. Our score is at most $H_w(x)$, because the penalty is nonnegative. Hence, $H_w(x)\ge T/2-C-m\log T\ge T/3$ for large $T$.

\emph{The counts are almost ordered.} Let $\zeta_w>0$ be the smallest positive difference between two coefficients of $w$, taking $\zeta_w=1$ if all coefficients are equal, and let $R_T=(2m\log T+4c_0)/\zeta_w$. We claim that
\begin{equation}\label{eq:almost-ordered}
 (x_k+1)\le (x_i+1)+R_T\qquad(i<k).
\end{equation}
If $w_i=w_k$ this is the order chosen within ties. Let $w_i>w_k$ and $x_k>x_i$. Swapping the two coefficients gives a direction of $\mathcal A$ with the same penalty function, whose height is at least $H_w(x)+g-2c_0$ with $g=(w_i-w_k)(x_k-x_i)\ge\zeta_w(x_k-x_i)$. By \eqref{eq:half-charge} its score exceeds ours by at least $g/2-2c_0$. This is at most $m\log T$, which gives \eqref{eq:almost-ordered}. In particular $T\le d((x_1+1)+R_T)$, so $T^{1-\varepsilon}\le T/(2d)\le (x_1+1)\le T$ for large $T$, and the first coordinate satisfies $|\log_T (x_1+1)-w_1|\le\varepsilon$.

\emph{Comparing with a changed suffix.} Let $j\ge2$ and $t=w_j$. We compare $w$ with competitors in $\mathcal A$ that have the same coefficients as $w$ up to some position $i\le j$ and different ones afterwards. If $t$ is too large for the count $x_j$, then replacing the coefficients from position $j$ onward by zeros saves more penalty than it loses in height. If $t$ is too small, then raising these coefficients gains more height than it costs in penalty. In both cases the score of the competitor would exceed the score of $w$ by more than $m\log T$, which is excluded by the assumption on $w$. To make this precise, let $P$ be the common charge function of the unchanged coefficients. Let $w'$ be such a competitor. Before accounting for either suffix, the score gain is
\begin{equation}\label{eq:shared-charge}
 \bigl(H_{w'}(x)-P(H_{w'}(x))\bigr)-\bigl(H_w(x)-P(H_w(x))\bigr)=\bigl(H_{w'}(x)-H_w(x)\bigr)-\bigl(P(H_{w'}(x))-P(H_w(x))\bigr).
\end{equation}
If height rises, at least half the gain survives by \eqref{eq:half-charge}. If it falls, $P(H_{w'}(x))\le P(H_w(x))$, so the loss is at most $H_w(x)-H_{w'}(x)$. To this difference add the suffix charges saved at $H_w(x)$ and subtract the new ones at $H_{w'}(x)$.

\begin{itemize}
\item \emph{Lower bound on $(x_j+1)$.} If $t=0$ there is nothing to prove, since $(x_j+1)\ge1$. For $t>0$, keep the first $j-1$ coefficients, choose the tag $0$ at position $j$ and finish with zeros. The new suffix is free, and we save at least the charge $H_w(x)^c/(2d)$ for $t$. By \eqref{eq:almost-ordered} the height loss is at most $\sum_{k\ge j}(x_k+1)+2c_0\le d((x_j+1)+R_T)+2c_0$. The competitor's score is at most ours plus $m\log T$, so
\[
 \frac{H_w(x)^c}{2d}\le d((x_j+1)+R_T)+2c_0+m\log T.
\]
Since $H_w(x)\ge T/3$ and $0<c\le1$, this reads $T^c\le C((x_j+1)+\log T+1)$, and therefore $(x_j+1)\ge T^{c-\varepsilon}$ for large $T$. Each left endpoint used in $c$ lies within $\varepsilon$ of a tag at least $t$, so $c>t-\varepsilon$ and $\log_T (x_j+1)>t-2\varepsilon$.

\item \emph{Upper bound on $(x_j+1)$.} Let $i\le j$ be the first position with $w_i=t$. If $i=1$ then $t=1$ and $(x_j+1)\le T$. Otherwise $t<\bar t=w_{i-1}$. Take the next larger tag $t^+$ in the partition of $[0,\bar t]$, keep the first $i-1$ coefficients and repeat $t^+$ from position $i$ to the end. Every repetition is available as the last tag of the next partition. The interval for $t^+$ starts at $U$, the right endpoint for $t$. All new suffix exponents are therefore at most $U$. There are fewer than $d$ of them, and the competitor's height is $H_{w'}(x)\le T+c_0\le2T$, so their total cost is at most $H_{w'}(x)^U/2\le T^U$. No changed coefficient decreases, and the one at position $j$ rises by $t^+-t$, giving a height gain of at least $(t^+-t)(x_j+1)-2c_0$. Half this gain, less $2c_0$ and the new charges, is at most $m\log T$:
\[
 (t^+-t)(x_j+1)\le4c_0+2T^U+2m\log T\le(t^+-t)T^{U+\varepsilon}
\]
for large $T$. Hence, $\log_T (x_j+1)\le U+\varepsilon<t+2\varepsilon$.
\end{itemize}
Altogether $|\log_T (x_i+1)-w_i|<2\varepsilon\le E_w$ for every $i$, so the direction fits.
\end{proof}

\section{One function for every state}\label{sec:assembly}

We sum the functions $W_w$ over the selected directions:
\[
 V(x)=\sum_{w\in\mathcal A} W_w(x),\qquad W_w=e^{G_w}(1+D_w).
\]
After a jump we evaluate each summand at the new state with its own fixed direction and penalty function. By linearity, $\L V=\sum_w \L W_w$. Since $D_w$ is bounded, each summand is between fixed positive multiples of $e^{G_w}$. The all-ones direction belongs to $\mathcal A$, and its score is $T-o(T)$, because its correction is bounded (Theorem~\ref{thm:height}(i)) and its penalty is $o(h)$ (Proposition~\ref{prop:selection}). Hence, $V$ has finite sublevel sets. It remains to check that the negative drift of a highest-scoring direction pays for all possible gains from the directions that fail to fit.

\begin{proposition}[Drift of the global weight]\label{prop:global-drift}
Outside a finite subset of $\Gamma$, $\L V(x)\le-1$.
\end{proposition}
\begin{proof}
Write $G_*$ for the largest score and let $\nu=\max_{w\in\mathcal A}\nu_w$, where $\nu_w$ is the exponent in Lemma~\ref{lem:uniform-preparation}. For large $T$ every fitting direction has negative drift by Lemma~\ref{lem:uniform-preparation}, and a maximizing direction fits by Proposition~\ref{prop:selection} and supplies a decrease of at least $c e^{G_*}T^{-\nu}$. By Proposition~\ref{prop:selection}, a direction that fails to fit lies more than $(\nu+M+1)\log T$ below the maximum. Its term in $V$ is at most $C e^{G_*}T^{-\nu-M-1}$. Each jump increases it by at most a constant times its present value, by \eqref{eq:one-jump}. The total exit rate $q(x)=\sum_{y\to y'\in\Re}\lambda_{y\to y'}(x)$ is at most $CT^M$, since $\lambda_y(x)\le T^M$ and $\kappa_{y\to y'}(x)\le\kappa_{\max}$, so its possible gain per unit time is at most $C e^{G_*}T^{-\nu-1}$. Summing over the finitely many directions yields
\[
 \L V(x)\le e^{G_*}\bigl(-cT^{-\nu}+C T^{-\nu-1}\bigr)\le-1
\]
for large $T$, since $G_*\ge T-o(T)$. There are only finitely many states below any fixed $T$.
\end{proof}

\Needspace{7\baselineskip}
\paragraph{From negative drift to finite mean returns.}
We use the continuous-time Foster criterion of Hajek~\cite[\S6.9.2, Proposition~6.16(b)]{Hajek2015}.
\begin{proposition}[Continuous-time Foster criterion]\label{prop:foster}
For an irreducible chain on a countable state space with positive finite exit rates, suppose $V\ge0$ is finite, has finite sublevel sets, and satisfies $\L V\le-1$ outside a finite set, with $\L V$ finite on that set. Then the mean return time to every state after departure is finite.
\end{proposition}
In Hajek's notation the drift condition is $\L V\le-1+\theta\mathbf1_C$. Take $C$ to be the finite exceptional set and choose $\theta>0$ large enough that $1+\L V\le \theta$ on $C$. Finite sublevel sets allow us to apply the criterion before proving non-explosion.

\begin{proof}[Conclusion of the proof of Theorem~\ref{thm:recurrence}]
Our class is infinite and irreducible, so every state has an active reaction. Only finitely many reactions are possible at a state, each with finite rate. Thus, its exit rate is positive and finite, and the sum defining $\L V$ is finite there. The all-ones direction gave finite sublevel sets, and Proposition~\ref{prop:global-drift} gives the negative drift. Proposition~\ref{prop:foster} now gives finite mean returns. Starting at any $z\in\Gamma$, the chain therefore returns to $z$ almost surely after every departure. By the strong Markov property, the successive holds at $z$ are independent exponentials of the same rate $q(z)>0$. Their durations sum to infinity almost surely, ruling out explosion.
\end{proof}

\section*{Acknowledgments} 

We thank the OpenRefinery STEM Research Program for an AI tool grant that covered the ChatGPT subscription used in this work.

\appendix
\section{Notation and dependencies}\label{app:notation}

Fix the network, its rates and a class $\Gamma$. For each direction, Section~\ref{sec:height} proves that $s_w$ takes finitely many values for every $\beta_w>0$, then chooses $\beta_w$ and $R_w$. In Section~\ref{sec:races} the bound $K_w$ on the number of cost gaps above one fixes the discount exponent $\sigma_w$, which in turn fixes the accuracy $E_w$. The correction $D_w$ works for every penalty satisfying \eqref{eq:half-charge}. This is why Section~\ref{sec:selection} can choose the directions and penalties last.

\begingroup
\small
\renewcommand{\arraystretch}{1.15}
\begin{longtable}{@{}p{0.25\textwidth}p{0.54\textwidth}l@{}}
\toprule
Notation & Meaning & See\\
\midrule
\endfirsthead
\toprule
Notation & Meaning & See\\
\midrule
\endhead
$x$, $\Gamma$, $\L $ & State, communicating class and generator. & \S\ref{sec:model}, \S\ref{sec:main}\\
$y\to y'$, $\kappa_{y\to y'}(x)$, $\lambda_y(x)$ & Reaction, rate coefficient, and mass-action factor. & Defs.~\ref{def:network}--\ref{def:kinetics}\\
$q(x)$, $M$ & Total exit rate and maximum source size. & \S\ref{sec:local}, \S\ref{sec:assembly}\\
$(n,k)$ & Counts of $A$ and $B$ in the running example. & Ex.~\ref{ex:running}\\
$w$, $\dtp{w,y}$ & Direction and order of a complex. & Def.~\ref{def:slope}\\
$x-y$, $\dtp{w,x-y}$ & Unused stock and unused level of a step $(x,y\to y')$. & Def.~\ref{def:step}\\
$\beta_w$, $R_w$ & Additive constant in the terminal charge $\dtp{w,z}+\beta_w$, and height threshold for a preparation. & Lem.~\ref{lem:terminal-choice}\\
$H_w(x)$, $s_w(x)=\dtp{w,x}-H_w(x)$ & Height and bounded state correction. & Def.~\ref{def:height}\\
$\dtp{w,y}-s_w(x)$ & Relative order of a step $(x,y\to y')$, equal to $H_w(x)-\dtp{w,x-y}$. & Def.~\ref{def:order}\\
$\gamma_w$ & Lower bound on the absolute value of the negative relative orders and on the height drop in an order-zero decrease. & Thm.~\ref{thm:height}\\
Preparatory step & Preserves height and has nonnegative relative order. & Def.~\ref{def:retained}\\
Preparation & Preparatory steps to a state with an active order-zero decrease. & Def.~\ref{def:preparation}\\
$\pathcost_w(x)$ & Least cost of a preparation from $x$, including one for the decrease that follows. & Def.~\ref{def:cost}\\
$\ell_w$, $K_w$ & Bounds on the preparation length needed for a minimum and on the number of minimum costs above one at a fixed height. & Lem.~\ref{lem:route-distance}\\
$\sigma_w$, $\delta_w$, $D_w(x)$ & Exponent discounting successive cost gaps, constant multiplier and resulting preparation correction. & \eqref{eq:preparation-correction}\\
$T=\sum_i(x_i+1)$ & Total size of the shifted state $x+1$. & \S\ref{sec:main}\\
$E_w$ & Allowed error in fitting $w$ to $\log_T (x+1)$. & Def.~\ref{def:fit}\\
$c_w$, $\nu_w$ & Constant and exponent in the negative drift estimate. & \eqref{eq:preparation-drift}\\
$\mathcal A$, $L,U,t$ & Finite direction family; endpoints and tag of a coefficient interval. & Prop.~\ref{prop:selection}\\
$b_w(h)$, $G_w(x)$ & Penalty and score $H_w(x)-b_w(H_w(x))$. & Def.~\ref{def:score}\\
$P$ & Charge function of the coefficients shared by a direction and its competitor. & \eqref{eq:shared-charge}\\
$W_w$, $V$ & Function $e^{G_w}(1+D_w)$ and sum over the selected directions. & \S\ref{sec:assembly}\\
\bottomrule
\end{longtable}
\endgroup


\begin{thebibliography}{10}

\bibitem{Agazzi2025}
A. Agazzi, D. F. Anderson, D. Cappelletti, L. Laurence, and J. C. Mattingly.
A proof of the chemical recurrence conjecture in two dimensions.
\href{https://arxiv.org/abs/2611.31221v1}{arXiv:2611.31221v1}, 23 November 2026.

\bibitem{ADE2018}
A. Agazzi, A. Dembo, and J.-P. Eckmann.
Large deviations theory for Markov jump models of chemical reaction networks.
\emph{Annals of Applied Probability} \textbf{28}(3), 1821--1855, 2018.
\href{https://doi.org/10.1214/17-AAP1344}{doi:10.1214/17-AAP1344}.

\bibitem{ADEgeometry2018}
A. Agazzi, A. Dembo, and J.-P. Eckmann.
On the geometry of chemical reaction networks: Lyapunov function and large deviations.
\emph{Journal of Statistical Physics} \textbf{172}(2), 321--352, 2018.
\href{https://doi.org/10.1007/s10955-018-2035-8}{doi:10.1007/s10955-018-2035-8}; \href{https://arxiv.org/abs/1710.07240v2}{arXiv:1710.07240v2}.

\bibitem{AgazziLaurence2026}
A. Agazzi and L. Laurence.
Bulk stability of bimolecular stochastic chemical reaction networks.
\href{https://arxiv.org/abs/2610.31550v1}{arXiv:2610.31550v1}, 22 October 2026.

\bibitem{ACKK2018}
D. F. Anderson, D. Cappelletti, M. Koyama, and T. G. Kurtz.
Non-explosivity of stochastically modeled reaction networks that are complex balanced.
\emph{Bulletin of Mathematical Biology} \textbf{80}, 2561--2579, 2018.
\href{https://doi.org/10.1007/s11538-018-0473-8}{doi:10.1007/s11538-018-0473-8}.

\bibitem{ACK2020}
D. F. Anderson, D. Cappelletti, and J. Kim.
Stochastically modeled weakly reversible reaction networks with a single linkage class.
\emph{Journal of Applied Probability} \textbf{57}(3), 792--810, 2020.
\href{https://doi.org/10.1017/jpr.2020.28}{doi:10.1017/jpr.2020.28}.

\bibitem{ACKN2020}
D. F. Anderson, D. Cappelletti, J. Kim, and T. D. Nguyen.
Tier structure of strongly endotactic reaction networks.
\emph{Stochastic Processes and their Applications} \textbf{130}(12), 7218--7259, 2020.
\href{https://doi.org/10.1016/j.spa.2020.07.012}{doi:10.1016/j.spa.2020.07.012}.

\bibitem{ACK2010}
D. F. Anderson, G. Craciun, and T. G. Kurtz.
Product-form stationary distributions for deficiency zero chemical reaction networks.
\emph{Bulletin of Mathematical Biology} \textbf{72}, 1947--1970, 2010.
\href{https://doi.org/10.1007/s11538-010-9517-4}{doi:10.1007/s11538-010-9517-4}.

\bibitem{AK2018}
D. F. Anderson and J. Kim.
Some network conditions for positive recurrence of stochastically modeled reaction networks.
\emph{SIAM Journal on Applied Mathematics} \textbf{78}(5), 2692--2713, 2018.
\href{https://doi.org/10.1137/17M1161427}{doi:10.1137/17M1161427}.

\bibitem{AndersonKurtz2015}
D. F. Anderson and T. G. Kurtz.
\emph{Stochastic Analysis of Biochemical Systems}.
Mathematical Biosciences Institute Lecture Series, Stochastics in Biological Systems, vol.~1.2. Springer, 2015.
\href{https://doi.org/10.1007/978-3-319-16895-1}{doi:10.1007/978-3-319-16895-1}.

\bibitem{CNP2013}
G. Craciun, F. Nazarov, and C. Pantea.
Persistence and permanence of mass-action and power-law dynamical systems.
\emph{SIAM Journal on Applied Mathematics} \textbf{73}(1), 305--329, 2013.
\href{https://doi.org/10.1137/100812355}{doi:10.1137/100812355}.

\bibitem{Dickson1913}
L. E. Dickson.
Finiteness of the odd perfect and primitive abundant numbers with $n$ distinct prime factors.
\emph{American Journal of Mathematics} \textbf{35}(4), 413--422, 1913.
\href{https://doi.org/10.2307/2370405}{doi:10.2307/2370405}.

\bibitem{GMS2014}
M. Gopalkrishnan, E. Miller, and A. Shiu.
A geometric approach to the global attractor conjecture.
\emph{SIAM Journal on Applied Dynamical Systems} \textbf{13}(2), 758--797, 2014.
\href{https://doi.org/10.1137/130928170}{doi:10.1137/130928170}.

\bibitem{Hajek2015}
B. Hajek.
\emph{Random Processes for Engineers}.
Cambridge University Press, 2015.
\href{https://hajek.ece.illinois.edu/Papers/randomprocJuly14.pdf\#page=215}{Proposition~6.16(b), p.~205 of the author's preproduction copy}.

\bibitem{Kriebel2026} A. Kriebel. Positive recurrence for single-linkage bimolecular weakly reversible stochastic reaction networks. Preprint, version 1.2.4, 2026. \href{https://doi.org/10.5281/zenodo.22089551}{doi:10.5281/zenodo.22089551}.

\bibitem{Norris1997}
J. R. Norris.
\emph{Markov Chains}.
Cambridge Series in Statistical and Probabilistic Mathematics. Cambridge University Press, 1997.

\bibitem{OpenAIBP2026}
OpenAI.
\emph{Boundedness and persistence of weakly reversible mass-action systems}.
Preprint, 25 September 2026, \href{https://github.com/openai/math/blob/adc7f1241b42e322a6451854ab7e4b4c146bf78a/preprints/Boundedness-and-persistence-of-weakly-reversible-mass-action-systems-September-25-2026/build/sections/02-affine.tex}{Lemma~2.1 and Corollary~2.2, Affine labels with a prescribed tolerance}.
Source repository at \url{https://github.com/openai/math}, accessed 23 November 2026.

\bibitem{OpenAI2026}
OpenAI.
\emph{Uniform permanence in weakly reversible mass-action systems}.
Preprint, 5 October 2026.
Source repository at \url{https://github.com/openai/math}, accessed 10 November 2026.

\bibitem{PCK2014}
L. Paulev\'e, G. Craciun, and H. Koeppl.
Dynamical properties of discrete reaction networks.
\emph{Journal of Mathematical Biology} \textbf{69}(1), 55--72, 2014.
\href{https://doi.org/10.1007/s00285-013-0686-2}{doi:10.1007/s00285-013-0686-2}.

\bibitem{Toth2018}
J. T\'oth, A. L. Nagy, and D. Papp.
\emph{Reaction Kinetics: Exercises, Programs and Theorems. Mathematica for Deterministic and Stochastic Kinetics}.
Springer, New York, 2018.
\href{https://doi.org/10.1007/978-1-4939-8643-9}{doi:10.1007/978-1-4939-8643-9}.

\bibitem{WX2023}
C. Wiuf and C. Xu.
Classification and threshold dynamics of stochastic reaction networks.
\href{https://arxiv.org/abs/2012.07954v3}{arXiv:2012.07954v3}, revised 24 January 2023.

\bibitem{Xu2026regularity}
C. Xu.
On the regularity of reaction systems.
\href{https://arxiv.org/abs/2409.05340v3}{arXiv:2409.05340v3}, revised 12 October 2026.

\bibitem{Xu2026firstorder}
C. Xu.
Exponential ergodicity of first order endotactic stochastic reaction systems.
\href{https://arxiv.org/abs/2601.00176v1}{arXiv:2601.00176v1}, 1 January 2026.

\end{thebibliography}
\end{document}